\documentclass[11pt]{article}

\title{\sc An explicit minorant for the amenability constant of the Fourier algebra}
\author{\sc Y. Choi}
\date{2nd December 2022}

\usepackage[a4paper,left=27mm,right=27mm,top=30mm,bottom=30mm]{geometry}
\usepackage[colorlinks=true,citecolor=red,linkcolor=blue]{hyperref}
\usepackage{sectsty}
\allsectionsfont{\sffamily}

\usepackage{amsmath,amssymb,amsfonts}
\usepackage[dvipsnames]{xcolor}

\numberwithin{equation}{section}
\usepackage{amsthm}
\newcounter{pulse}[section]
\numberwithin{pulse}{section}

\newcommand{\newheadfont}{\bfseries} 

\newtheoremstyle{newplain} 
  {\topsep}   
  {\topsep}   
  {\itshape}  
  {0pt}       
  {\newheadfont} 
  {.}         
  {5pt plus 1pt minus 1pt} 
  {}          

\newtheoremstyle{newdef} 
  {\topsep}   
  {\topsep}   
  {\normalfont}  
  {0pt}       
  {\newheadfont} 
  {.}         
  {5pt plus 1pt minus 1pt} 
  {}          

\newtheoremstyle{newrem}
  {\topsep}   
  {\topsep}   
  {\normalfont}  
  {0pt}       
  {\newheadfont} 
  {.}         
  {5pt plus 1pt minus 1pt} 
  {}          

\theoremstyle{newplain}
\newtheorem{thm}[pulse]{Theorem}
\newtheorem{prop}[pulse]{Proposition}
\newtheorem{lem}[pulse]{Lemma}
\newtheorem{cor}[pulse]{Corollary}
\theoremstyle{newdef}
\newtheorem{dfn}[pulse]{Definition}

\newtheorem{eg}[pulse]{Example}
\theoremstyle{newrem}
\newtheorem{rem}[pulse]{Remark}

\newtheorem{qn}[pulse]{Question}

\newenvironment{romnum}{%
\begin{enumerate}

}{\end{enumerate}\ignorespacesafterend}

\newcommand{\defeq}{:=}
\newcommand{\dt}[1]{{\sf {#1}}}

\newcommand{\Cplx}{\mathbb C}
\newcommand{\Nat}{\mathbb N}

\newcommand{\Real}{\mathbb R}
\newcommand{\Zahl}{\mathbb Z}

\newcommand{\bbT}{\mathbb T}

\newcommand{\AM}{\operatorname{AM}} 
\newcommand{\AD}{\operatorname{AD}} 

\newcommand{\muad}{\mu^{\rm(ad)}}
\newcommand{\nud}{\nu^{\rm(d)}}

\newcommand{\Bdd}{{\mathcal B}}
\newcommand{\Tr}{\operatorname{Tr}} 

\newcommand{\UT}[3]{\begin{pmatrix} 1 & #1 & #3 \\ 0 & 1 & #2 \\ 0 & 0 & 1 \end{pmatrix}}

\newcommand{\veps}{\varepsilon}

\newcommand{\gm}{\gamma}
\newcommand{\lm}{\lambda}

\newcommand{\Gm}{\Gamma}
\newcommand{\Lm}{\Lambda}
\newcommand{\Om}{\Omega}

\newcommand{\Heis}{{\mathbb H}}

\newcommand{\cH}{{\mathcal H}}

\newcommand{\cU}{{\mathcal U}}

\newcommand{\sX}{{\mathsf X}}

\newcommand{\Ahat}{\widehat{A}}
\newcommand{\Bhat}{\widehat{B}}
\newcommand{\Dhat}{\widehat{D}}
\newcommand{\Ghat}{\widehat{G}}

\newcommand{\Delhat}{\widehat{\Delta}}
\newcommand{\Gmhat}{\widehat{\Gamma}}
\newcommand{\Lmhat}{\widehat{\Lambda}}

\newcommand{\twice}[1]{{#1\times#1}}
\newcommand{\longhat}[1]{\left(#1\right)^\wedge}
\newcommand{\hatofsquare}[1]{\longhat{\twice{#1}}}
\newcommand{\hatonehattwo}[1]{\widehat{{#1}_1}\times\widehat{{#1}_2}}
\newcommand{\hatonetwo}[1]{\longhat{{#1}_1\times{#1}_2}}

\newcommand{\ip}[2]{{\langle#1,#2\rangle}}

\newcommand{\wstar}{\ensuremath{{\rm w}^*}}

\newcommand{\tp}{\mathbin{\otimes}}
\newcommand{\ptp}{\mathbin{\widehat{\otimes}}}
\newcommand{\maxtp}{\mathbin{\otimes_{\max}}}

\newcommand{\Abs}[1]{{\left\vert#1\right\vert}}
\newcommand{\norm}[1]{{\Vert#1\Vert}}

\newcommand{\FA}{\operatorname{A}}
\newcommand{\Cst}{\operatorname{C}^*}
\newcommand{\FS}{\operatorname{B}}

\newcommand{\fsnorm}[2][]{{\norm{#2}_{{\rm B}#1}}}   

\newcommand{\diag}{\operatorname{diag}} 

\newcommand{\adiag}{\operatorname{adiag}}

\newcommand{\maxdeg}{\operatorname{maxdeg}}

\renewcommand{\Re}{\operatorname{Re}}

\newcommand{\comm}{\operatorname{comm}} 

\newcommand{\twomat}[4]{\begin{pmatrix} #1 & #2 \\ #3 & #4 \end{pmatrix}}
\newcommand{\grindex}[2]{{|#1\colon #2|}}
\newcommand{\Zindex}[1]{\grindex{#1}{Z(#1)}}

\usepackage{tikz-cd}

\renewcommand{\newheadfont}{\sc}  

\begin{document}

\maketitle

\begin{abstract}
We show that if a locally compact group $G$ is non-abelian then the amenability constant of its Fourier algebra is $\geq 3/2$, extending a result of Johnson (JLMS, 1994) who proved that this holds for finite non-abelian groups. Our lower bound, which is known to be best possible, improves on results by previous authors and answers a question raised by Runde (PAMS, 2006).

To do this we study a minorant for the amenability constant, related to the anti-diagonal in $G\times G$, which was implicitly used in Runde's work but hitherto not studied in depth. Our main novelty is an explicit formula for this minorant when $G$ is a countable virtually abelian group, in terms of the Plancherel measure for $G$.
 As further applications, we characterize those non-abelian groups where the minorant attains its minimal value, and present some examples to support the conjecture that  the minorant always coincides with the amenability constant.



\medskip\noindent
MSC 2020:
46H20
(primary);
43A30,
46J10
(secondary)
\end{abstract}

%
%

\bigskip
\paragraph{Note added for arXiv v5.}
While reading the page proofs for the final published version, I~noticed some remaining typos (mostly errors of punctuation). I have fixed these errors here. On the other hand, this version does not incorporate all of the copy-editors' changes, some of which were ``improvements'' to the grammar in my original wording. So this version does not replace the official version of record. Nevertheless, the mathematical substance should be identical to the published version.

\hfill --- YC, Lancaster, 28th June 2023

 \bigskip\hrule\newpage
%

\begin{section}{Introduction}\label{s:introduction}

\begin{subsection}{Background context}\label{ss:context}
The Fourier algebra of a locally compact group $G$, denoted by $\FA(G)$, is a Banach algebra of functions on $G$ whose norm encodes the group structure of $G$. By a celebrated theorem of Walter \cite{Wal_JFA72},
$\FA(G)$ is a complete invariant for $G$ in the following sense: given two locally compact groups $G_1$ and $G_2$ the Banach algebras $\FA(G_1)$ and $\FA(G_2)$ are isometrically algebra-isomorphic if and only if $G_1$ and $G_2$ are topologically group-isomorphic. Thus $\FA(G)$ is a finer invariant than the commutative $\Cst$-algebra $C_0(G)$, which only remembers the underlying topological space of~$G$ (and in particular cannot distinguish between discrete groups of the same cardinality).

It is therefore natural to study various invariants of Banach algebras, specialized to the setting of Fourier algebras, to see what information they remember about the group one starts with. Typically, algebraically motivated invariants of Banach algebras are \emph{isomorphic} rather than \emph{isometric} invariants, and hence cannot distinguish between groups whose Fourier algebras are isomorphic as topological algebras. For instance, if $G_1$ and $G_2$ are finite groups with the same cardinality~$n$, then $\FA(G_1)$ and $\FA(G_2)$ are both (non-isometrically) isomorphic as unital Banach algebras to $(\Cplx^n, \norm{\quad}_\infty)$. In particular, since the latter Banach algebra is known to be \dt{amenable}, so are $\FA(G_1)$ and $\FA(G_2)$, and the notion of amenability cannot distinguish between $\FA(G_1)$ and $\FA(G_2)$ even though $G_1$ and $G_2$ may have very different group structure.

However, if $A$ is an amenable Banach algebra one may consider its \dt{amenability constant} $\AM(A)$. This is an isometric invariant, and so has the potential to distinguish between two Fourier algebras that happen to be non-isometrically isomorphic; moreover, it is a numerical invariant, and hence lends itself to estimates and inequalities once one exploits its functorial properties.
Even for finite groups, it displays interesting behaviour, as demonstrated by Johnson in his study \cite{BEJ_AG} of the (non-)amenability properties of Fourier algebras of \emph{compact} groups.

Specifically, when $G$ is a finite group, \cite[Theorem 4.1]{BEJ_AG} provides an explicit and unexpectedly simple formula for $\AM(\FA(G))$:
\begin{equation}\label{eq:BEJ formula}
\AM(\FA(G))= \frac{1}{|G|}\sum_{\pi\in\Ghat} (d_\pi)^3
\end{equation}
where $\Ghat$ denotes the set of (equivalence classes of) irreducible representations of $G$, and $d_\pi$ denotes the dimension of the representation ~$\pi$. As we will see in Theorem \ref{t:AD of finite}, the right-hand side of \eqref{eq:BEJ formula} has independent meaning, which seems to have been overlooked in \cite{BEJ_AG} and the work of subsequent authors.

The formula \eqref{eq:BEJ formula} has several applications. For instance, it shows that $\AM(\FA(G))=1$ for any \emph{finite abelian}~$G$. With the aid of some basic group theory it also implies the following ``gap theorem'' for the possible values of $\AM(\FA(G))$ in the class of finite groups.

\begin{thm}[{\cite[Proposition 4.3]{BEJ_AG}}]
\label{t:BEJ gap}
If $G$ is a finite non-abelian group then $\AM(\FA(G))\geq 3/2$. Equality can be achieved (for instance, take $G$ to be the dihedral group of order $8$ or the quaternion group).
\end{thm}

Another consequence of \eqref{eq:BEJ formula}, already observed by Johnson \cite[Corollary 4.2]{BEJ_AG}, is the following formula:
\begin{equation}\label{eq:AMAG-product}
\AM(\FA(G_1\times G_2))=\AM(\FA(G_1)) \AM(\FA(G_2)) \quad\text{for any finite groups $G_1$ and $G_2$.}
\end{equation}
This formula was used in conjunction with Theorem \ref{t:BEJ gap} to establish the non-amenability of the Fourier algebras of a large class of profinite groups (see \cite[Theorem 4.5]{BEJ_AG}).

\begin{rem}\label{r:AMAG-of-product}
Equation \eqref{eq:AMAG-product} is not obvious from the general definition of the amenability constant of a Banach algebra. It is true that for arbitrary Banach algebras $A$ and $B$ one has $\AM(A\ptp B)\leq \AM(A)\AM(B)$, where $\ptp$ denotes the projective tensor product; and one can show that equality holds if $A$ and $B$ are commutative, semisimple and finite-dimensional. However, this would only tell us about $\FA(G_1)\ptp\FA(G_2)$, which in general is not isometrically isomorphic to $\FA(G_1\times G_2)$.
Indeed, it seems to be unknown whether the formula in Equation \eqref{eq:AMAG-product} remains true when $G_1$ and $G_2$ are general locally compact groups.
\end{rem}

There are difficulties in extending \cite[Theorem 4.1]{BEJ_AG} to the non-compact setting.
Johnson's proof of the explicit formula \eqref{eq:BEJ formula} is indirect, since it proceeds by analyzing an auxiliary Banach algebra $\FA_\gm(G)$ and showing that
for finite $G$ a particular idempotent in $\FA_\gm(G)$ is mapped isometrically to the unique (virtual) diagonal for $\FA(G)$. Although the definition of $\FA_\gm(G)$ can be extended to locally compact virtually abelian~$G$, the relationship with $\FA(G)$ is less clear in non-compact cases. Moreover, when $G$ is infinite, one no longer has a unique virtual diagonal for $\FA(G)$.

\end{subsection}

\begin{subsection}{A minorant for the amenability constant}
The study of (non-)amenability in \cite{BEJ_AG} relied heavily on the matrix-valued Fourier transform available for compact groups.
Using a different perspective, Forrest and Runde proved that $\FA(G)$ is amenable if and only if $G$ is \emph{virtually abelian} \cite[Theorem 2.3]{FR_MZ05}. A key step in their proof is the following result: 
if $\FA(G)$ is amenable, then the indicator function of the set $\adiag(G)\defeq\{ (x,x^{-1}) \colon x\in G\}$ belongs to the Fourier--Stieltjes algebra of $G_d\times G_d$.
Here $G_d$ denotes $G$ equipped with the discrete topology.

In fact, as observed by Runde in a follow-up paper \cite{Run_PAMS06}, the argument in \cite{FR_MZ05} yields a more precise statement: $\AM(\FA(G))$ is bounded below by the Fourier--Stieltjes norm of $1_{\adiag(G)}$.
This quantity, which we shall denote by $\AD(G_d)$,  is the ``minorant'' alluded to in the title; it can be defined without reference to amenability, and depends only on the underlying group of $G$ rather than its topological group structure.

\begin{rem}\label{r:AD-better-than-AMA}
As an invariant of $G$, $\AD(G_d)$ enjoys slightly better functorial properties than $\AM(\FA(G))$.
For instance, given two virtually abelian groups $\Gm_1$ and~$\Gm_2$, one has
$\AD(\Gm_1\times\Gm_2)= \AD(\Gm_1)\AD(\Gm_2)$; as previously remarked, the corresponding result for $\AM(\FA(\cdot))$ is only known in particular cases.
It also behaves well with respect to taking arbitrary subgroups, while $\AM(\FA(\cdot))$ is only known to behave well with respect to closed subgroups.
Proofs of these functorial properties will be given below (see Proposition \ref{p:AD gen prop}).
\end{rem}

Besides its intrinsic interest, $\AD$ is a natural tool if one seeks to obtain \emph{lower bounds} on $\AM(\FA(G))$ for particular $G$, which would be hard to ascertain directly from the definition of the amenability constant. This was already noted in \cite[Section~3]{Run_PAMS06}, where after some partial results the following question is raised.

\paragraph{Runde's question, paraphrased} (see also \cite[Remark 4.9]{Juselius_MSc}).
\
\begin{itemize}
\item (Strong form.)
Does the lower bound in Theorem \ref{t:BEJ gap} remain true for all locally compact non-abelian groups, not just the finite ones?
\item (Weak form.)
Can we find some $C>1$ such that $\AM(\FA(G))\geq C$ for every locally compact non-abelian group $G$?
\end{itemize}

Progress on these questions has been hampered by the lack of a formula such as \eqref{eq:BEJ formula} for infinite groups.
Nevertheless, the weak form of Runde's question was answered positively by Forrest and Runde, by using results on norms of idempotent Schur multipliers to obtain lower bounds on $\AD(G_d)$ (although they do not use this notation). In particular, it is observed in the remarks following \cite[Theorem 3.4]{FR_CMB11} that one can take $C=2/\sqrt{3}\approxeq 1.1547$. 
An improved value of $C=9/7\approxeq 1.2857$ can be extracted\footnotemark\ from results of Mudge and Pham \cite{MP_JFA16}; strictly speaking, this is not stated explicitly in their paper, but it follows from the proof of \cite[Theorem 2.2]{MP_JFA16} and observations before~it.
\footnotetext{The author is grateful to N. Juselius for sharing a copy of \cite{Juselius_MSc}, which provides a useful synopsis of these developments. The results of Forrest--Runde and Mudge--Pham actually prove something sharper, concerning lower bounds on cb-multiplier norms of functions rather than their Fourier--Stieltjes norms, but we shall not discuss this extra refinement here.}

Despite this progress, the original stronger version of Runde's question has remained unanswered until now. Naively, since $\AM(\FA(G))\geq \AM(\FA(H))$ whenever $H$ is a closed subgroup of~$G$, one might hope to prove that $\AM(\FA(G))\geq 3/2$ by locating non-abelian finite subgroups of $G$ and applying Theorem~\ref{t:BEJ gap}. However, as we shall see in Section \ref{ss:group background}, there are virtually abelian non-abelian groups in which every finite subgroup is abelian, and such groups cannot be handled with this strategy.

\end{subsection}

\begin{subsection}{Main results of this paper}
In this paper we develop new machinery which yields, as a by-product, a positive answer to the strong form of Runde's question.
By \cite[Theorem 2.3]{FR_MZ05}, it is enough to restrict our attention to virtually abelian groups, where more tools are available.
Rather than using the perspective of operator theory, and trying to find improved estimates of Schur multiplier norms, we use the theory of the operator-valued Fourier transform. (In a sense this is more in the spirit of Johnson's original paper~\cite{BEJ_AG}, although as indicated in Section \ref{ss:context} we have to go beyond the techniques of \cite{BEJ_AG} once we leave the realm of compact groups.)

If one wishes to check how tight $\AD(G_d)$ is as a lower bound for $\AM(\FA(G))$, then a natural place to start is with finite groups where Johnson's formula~\eqref{eq:BEJ formula} is available. For finite $G$, the Fourier--Stieltjes algebra of $G_d\times G_d$ is just the Fourier algebra of $G\times G$, and one can calculate norms in this algebra in terms of the matrix-valued Fourier transform for $G\times G$. Working through the details, it turns out that the Fourier coefficients of $1_{\adiag(G)}$ have a particular simple form, yielding a surprising result.

\begin{thm}[Johnson's formula, revisited]
\label{t:AD of finite}
Let $G$ be a finite group. Then
\begin{equation}\label{eq:AD of finite}
\AD(G) = \frac{1}{|G|} \sum_{\pi\in\Ghat} (d_\pi)^3 = \AM(\FA(G)). 
\end{equation}
\end{thm}

Our result shows that \cite[Theorem 4.1]{BEJ_AG} should not merely be regarded as a numerical coincidence, where $\AM(\FA(G))$ happens to be equal to some ``random number'' defined in terms of $\Ghat$, but rather that an inequality relating two {\it a priori} different invariants for arbitrary $G$ is actually an equality for all finite groups.

When $\Gm$ is infinite and virtually abelian, one still has a version of the Plancherel formula and hence a version of the Fourier inversion theorem. By using the explicit formula for the Fourier coefficients of $1_{\adiag(G)}$ in the finite case, and making a suitable {\itshape Ansatz},  we obtain an explicit formula for $\AD(\Gm)$, at least when $\Gm$ is countable.

\begin{thm}[The main formula]
\label{t:main formula}
Let $\Gm$ be countable and virtually abelian, with Plancherel measure $\nu$ normalized so that $1 = \int_{\Gmhat} d_\pi\,d\nu(\pi)$.
 For each $n\in\Nat$, let $\Om_n = \{ \pi\in \Gmhat \colon d_\pi = n\}$. Then
\begin{equation}\label{eq:AD-explicit-formula}
\AD(\Gm) = \int_{\Gmhat} (d_\pi)^2 \,d\nu(\pi) = \sum_{n\in\Nat} n^2 \nu(\Om_n).
\end{equation}
\end{thm}

Although  the new part of Theorem \ref{t:AD of finite} is a special case of Theorem \ref{t:main formula}, it seems more natural to state the results separately and prove the special case first, since this approach provides crucial motivation for the key {\itshape Ansatz} in the proof of the general result.
Moreover, the proof of the general result is technically more involved, relying on a non-abelian version of Bochner's theorem that was developed in work of Arsac \cite[Section 3G]{Ars_Lyon76}.

Combining the formula \eqref{eq:AD-explicit-formula} with an upper bound on $\nu(\Om_1)$, and exploiting a ``monotonicity principle'' for $\AD$, we are able to improve the existing lower bounds on $\AM(\FA(G))$ for locally compact non-abelian groups, and show that the ``strong form'' of Runde's question has a positive answer.

\begin{thm}[Sharp lower bound]
\label{t:sharp LB}
If $\Gm$ is a (discrete) non-abelian group then $\AD(\Gm)\geq 3/2$. Consequently $\AM(\FA(G))\geq 3/2$ for every locally compact non-abelian group~$G$.
\end{thm}

With more work, we can go further and determine exactly which non-abelian groups attain the lower bound on $\AD(\Gm)$ in Theorem \ref{t:sharp LB}. Note that the next result does not require any countability assumptions.

\begin{thm}[Non-abelian groups where $\AD$ attains its minimum]
\label{t:minimizers of AD}
Let $\Gm$ be a (discrete) group. Then $\AD(\Gm)=3/2$ if and only if $\Zindex{\Gm}=4$.
\end{thm}

\begin{rem}
If $\Zindex{\Gm}=4$ then basic group theory implies that $\Gm/Z(\Gm)$ is non-cyclic and hence is isomorphic to $C_2\times C_2$ (see Lemma \ref{l:was-prop-F} below). Thus the groups described by Theorem \ref{t:minimizers of AD} are exactly the non-trivial central extensions of $C_2\times C_2$. In principle they can be classified by standard methods from group cohomology, although we have not found an explicit list in the literature.
\end{rem}

\begin{rem}
We may combine Theorems \ref{t:AD of finite} and \ref{t:minimizers of AD} to obtain the following result, which has not appeared before in the literature:
\textit{given a finite group $G$, $\AM(\FA(G))=3/2$ if and only if $\Zindex{G}=4$.}
However, it is also possible to derive this result directly from Johnson's formula~\eqref{eq:AD of finite}, using some basic properties of character tables. Details of this character-theoretic argument will appear in future work, as part of a more systematic investigation of $\AM(\FA(G))$ for finite groups.
\end{rem}

In view of Theorem \ref{t:AD of finite}, it is natural to ask if the inequality $\AD(G_d)\leq\AM(\FA(G))$ is always an equality. While we are unable to resolve this question at present, we can show that equality holds for certain infinite families of semidirect product groups (see also Remark~\ref{r:KFT_crystal}).

To put these examples in context, we return briefly to the results of \cite{Run_PAMS06}. If $G$ is a locally compact virtually abelian group, then $\maxdeg(G)\defeq \sup\{ d_\pi \colon \pi\in \Ghat\}$ is finite. Runde showed that for such $G$ we have a chain of inequalities
\begin{equation}\label{eq:AD-AM-maxdeg}
\AD(G_d)\leq \AM(\FA(G))\leq \maxdeg(G).
\end{equation}
(The first inequality is proved in \cite[Lemma 3.1]{Run_PAMS06} and the second inequality is proved in \cite[Theorem 2.7]{Run_PAMS06}.)
We can now state our final main result.

\begin{thm}[Sufficient conditions for equality to hold in \eqref{eq:AD-AM-maxdeg}]
\label{t:all equal}
Let $N$ be a torsion-free locally compact abelian group,
and let $p$ be a prime. Let $\alpha$ be a non-trivial action of the cyclic group~$C_p$ on $N$, and let $G=N\rtimes_\alpha C_p$. Then
\begin{equation}\label{eq:AD=AM=maxdeg}
\AD(G_d)=\AM(\FA(G))=\maxdeg(G)=p.
\end{equation}
\end{thm}

These are the first recorded examples of groups not of the form (abelian)$\times$(finite), where the amenability constant of the Fourier algebra is finite and can be calculated explicitly.

\begin{rem}\
\begin{romnum}
\item
The groups covered by Theorem \ref{t:all equal} include examples which are discrete, (non-connected) Lie, or compact: let $V$ be either $\Zahl$, $\Real$ or the $2$-adic integers respectively, and then take $N=V^p$ equipped with the natural permutation action of $C_p$.
\item
If $G$ is a finite non-abelian group, then it follows immediately from Johnson's formula \eqref{eq:BEJ formula} that $\AM(\FA(G))< \maxdeg(G)$. Hence the phenomenon in Theorem \ref{t:all equal} is special to the setting of infinite~$G$.
\end{romnum}
\end{rem}

\begin{rem}[Crystal groups]\label{r:KFT_crystal}
The proof of Theorem \ref{t:all equal} uses the fact that a suitably chosen countable subgroup $\Gamma \leq G$ has the following property: the set $\{ \pi \in \Gmhat\colon d_\pi < \maxdeg(\Gm)\}$ has Plancherel measure zero. In response to an earlier version of this paper, K. F. Taylor has pointed out (personal communication) that the same property holds for any \dt{crystal group}, i.e.\ a discrete subgroup $\Gm\ < {\rm Isom}(\Real^n)$ such that the orbit space $\Real^n/\Gm$ is compact. Using this fact, it follows easily from Theorem \ref{t:main formula} that for a crystal group $\Gm$ we also have $\AD(\Gm)=\AM(\FA(\Gm))=\maxdeg(\Gm)$. Moreover, $\maxdeg(\Gm)$ is equal to the index in $\Gm$ of its translation subgroup (equivalently, equal to the order of the \dt{point group} of $\Gm$).
\end{rem}

After some preliminaries in Section \ref{s:prelim}, the proofs of Theorems
\ref{t:AD of finite},
\ref{t:main formula},
\ref{t:sharp LB},
\ref{t:minimizers of AD}
and
\ref{t:all equal}
will be given in the following sections.
We finish the paper with some questions and suggestions for future work.
In an appendix, we present some results which relate the dual of a product of two locally compact groups with the product of the respective duals, not just at the level of sets but at the level of topological spaces; this is based on a suggestion by an anonymous referee.\footnote{See also the ``author's note'' at the end of Section \ref{ss:dual-of-product}.} These results work for general locally compact groups, without requiring extra conditions such as second-countability or being Type~I.

\end{subsection}

\end{section}

\begin{section}{Preliminaries}
\label{s:prelim}

\begin{subsection}{Group-theoretic background}\label{ss:group background}
A group $G$ is said to be \dt{virtually abelian} (VA for short) if it has an abelian subgroup of finite index.
For convenience, we shall also use the abbreviation VANA, which stands for ``virtually abelian non-abelian''.

Note that this property does not depend on any topological group structure that $G$ may or may not possess.
Note also that subgroups of VA groups are VA:
if $H$ has finite index in $G$ and $\Gm\leq G$ then $\grindex{\Gm}{H\cap\Gm} \leq \grindex{G}{H}$.

Clearly every finite group is VA, as is every abelian group. However the class of VA groups displays diverse behaviour, much more than the ``trivial'' examples of the form (abelian)$\times$(finite) might suggest. We present two examples, which are close to abelian in the sense that they have abelian subgroups of index $2$, and in which every non-trivial finite subgroup has order $2$ (so in particular, is abelian).

\begin{eg}[Generalized dihedral groups]\label{eg:dihedral-ish}
Let $N$ be any torsion-free LCA group and let $G=N\rtimes C_2$, where $C_2$ acts on $N$ by $x\mapsto -x$ (using additive notation). Regard $N$ as a normal subgroup of $G$, via the canonical embedding $N\to N\rtimes C_2$. Then every non-identity element of $N$ has infinite order in $G$, while every element of $G\setminus N$ is an involution; moreover, $xy\in N$ for all $x,y\in G\setminus N$.
Therefore the only finite subgroups of $G$ are the trivial subgroup $\{e\}$ and subgroups of the form $\{e, t\}$ for some $t\in G\setminus N$.
\end{eg}

The following example seems to be well known to specialists in group theory, and can be found in \cite[Section 6.8]{Folland_ed2}.
It may be viewed as a discrete analogue of the ``reduced Heisenberg group''.

\begin{eg}\label{eg:red-heis-mod2}
Consider the group
\[
\Heis = \left\{ \UT{a}{b}{c} \colon a,b,c\in\Zahl\right\} \leq {\rm SL}_3(\Zahl)
\]
whose elements will be denoted by triples $(a,b,c)$ for sake of conciseness.
Define $\Heis_2$ to be the quotient of $\Heis$ by the central subgroup $\{ (0,0,c)\colon c\in2\Zahl\}$.
We denote elements of $\Heis_2$ by triples $(a,b,[c])$ where $a,b\in\Zahl$ and $c\in\Zahl/2\Zahl$.
%
There is a surjective homomorphism $q: \Heis_2\to\Zahl^2$, given by $q(a,b,[c])= (a,b)$, and since ${\mathbb Z}^2$ is torsion-free, the only finite-order elements in $\Heis_2$ are those in $\ker(q)$.  Thus the only finite subgroups of $\Heis_2$ are the trivial subgroup and the 2-element subgroup $\ker(q)$. On the other hand, one can show by direct calculation that
\[ N\defeq\{(a,b,[c])\colon a\in \Zahl,b\in 2\Zahl, [c]\in \Zahl/2\Zahl\} \]
is an abelian subgroup of $\Heis_2$, and clearly $\grindex{\Heis_2}{N}=2$.
\end{eg}

In fact, for both of these examples, one can calculate the exact value of the anti-diagonal constant using the new results described in the introduction. The group $G$ from Example \ref{eg:dihedral-ish} falls into the class of groups covered by Theorem \ref{t:all equal}, and therefore satisfies $\AD(G)=2$. The group $\Heis_2$ in Example \ref{eg:red-heis-mod2} turns out to have centre $\{(a,b,[c])\colon a\in 2\Zahl,b\in 2\Zahl, [c]\in \Zahl/2\Zahl\}$, and quotienting out by the centre gives a copy of $C_2\times C_2$; so $\Heis_2$ falls into the class of groups covered by Theorem \ref{t:minimizers of AD}, and hence\footnotemark\ $\AD(\Heis_2)=3/2$.
\footnotetext{One can also give a more direct calculation of $\AD(\Heis_2)$ using the main formula in Theorem \ref{t:main formula}, see Example \ref{eg:AD red-heis-mod2} below.}

\end{subsection}

\begin{subsection}{Fourier and Fourier--Stieltjes algebras}\label{ss:AG-BG}
Here and throughout the paper: if $X$ is a set then $c_{00}(X)$ denotes the set of all finitely supported functions $X\to \Cplx$.

For the definition of the Fourier algebra $\FA(G)$ we refer the reader to \cite{eymard} or the recent text \cite{KL_AGBGbook}. We shall not discuss amenability of Banach algebras in this paper, as we are only concerned with the minorant  $\AD(G_d)$ rather than $\AM(\FA(G))$ itself. However, to study $\AD(G_d)$ we must first discuss norms of elements inside Fourier--Stieltjes algebras of discrete groups.
Since there are various equivalent definitions of the Fourier--Stieltjes algebra and equivalent definitions of its canonical norm, we review some terminology and notation.

An important convention throughout this paper is that the word ``representation'' \emph{always means} ``unitary representation''; this short-hand is merely to avoid unnecessary repetition.

\begin{dfn}
Let $\Gm$ be a discrete group. For a representation $\pi:\Gm\to \cU(H)$, a \dt{coefficient function of~$\pi$} is a function $\Gm\to\Cplx$ of the form $x\mapsto \ip{\pi(x)\xi}{\eta}$ where $\xi,\eta\in H$.
The \dt{Fourier--Stieltjes algebra of $\Gm$} is
\[
\FS(\Gm) \defeq \bigcup_\pi \{ \text{coefficient functions of $\pi$} \} \subseteq \ell^\infty(\Gm).
\]
where the union is taken over all representations of $\Gm$. Clearly $\FS(\Gm)$ is closed under scalar multiplication; by taking direct sums and tensor products of representations, we see that it is also closed under pointwise addition and multiplication.
\end{dfn}

Let $\Gm$ be a discrete group and let $b\in c_{00}(\Gm)$. If $\theta:\Gm\to\cU(H)$ is a representation, we write $\theta(b)$ for the (finite) sum $\sum_{x\in\Gm} b(x)\theta(x)\in\Bdd(H)$. 
We then define
$\norm{b}_{\Cst(\Gm)} \defeq \sup_\theta \norm{\theta(b)}$ where the supremum is taken over all representations of $\Gm$.

\begin{dfn}\label{d:fsnorm as sup}
For a discrete group $\Gm$ and $f\in \FS(\Gm)$, we define
\begin{equation}
\fsnorm{f}
\defeq \sup \left\{ \sum_{x\in\Gm} f(x)b(x) \colon b\in c_{00}(\Gm), \norm{b}_{\Cst(\Gm)}\leq 1 \right\}
\end{equation}
\end{dfn}

The following result is standard but included here for convenient reference.
\begin{lem}\label{l:fsnorm-as-inf}
Let $f\in \FS(\Gm)$. If $f$ can be expressed as a coefficient function $\ip{\theta(\cdot)\xi}{\eta}$ for some representation $\theta:\Gm\to\cU(H)$ and vectors $\xi,\eta\in H$, then $\fsnorm{f}\leq \norm{\xi}\norm{\eta}$.
Moreover, one can always find some choice of $H$, $\theta$, $\xi$ and $\eta$ for which  $f=\ip{\theta(\cdot)\xi}{\eta}$ and $\fsnorm{f}=\norm{\xi}\norm{\eta}$.
\end{lem}

The first part of the lemma follows from Definition \ref{d:fsnorm as sup} and the Cauchy--Schwarz inequality.
The second part follows by considering a suitable GNS representation and using polar decomposition; see e.g.\ \cite[Lemma 2.1.9]{KL_AGBGbook}.

\begin{lem}\label{l:FS basics}
\begin{romnum}
\item\label{li:restriction}
Let $\Lm$ be a discrete group and let $f\in \FS(\Lm)$.
If $\Lm_0$ is a subgroup of $\Lm$, then ${f\vert}_{\Lm_0}\in\FS(\Lm_1)$ and $\fsnorm{{f\vert}_{\Lm_0}}\leq \fsnorm{f}$.
\item\label{li:norm-of-tensor}
Let $\Lm_1$ and $\Lm_2$ be discrete groups, and let $f_i\in \FS(\Lm_i)$ for $i=1,2$. Define the function $f_1\tp f_2:\Lm_1\to\Lm_2\to \Cplx$ by $(f_1\tp f_2)(x,y)=f_1(x)f_2(y)$. Then $f_1\tp f_2\in \FS(\Lm_1\times\Lm_2)$ and
\[
\fsnorm{f_1\tp f_2} = \fsnorm{f_1} \fsnorm{f_2} \;.
\]

\end{romnum}
\end{lem}

\begin{proof}
Part \ref{li:restriction} is a straightforward consequence of Lemma \ref{l:fsnorm-as-inf}, which we leave to the reader.
The proof of part \ref{li:norm-of-tensor} requires more work.
By Lemma \ref{l:fsnorm-as-inf}, for $i=1,2$
we can write
$f_i=\ip{\theta_i(\cdot)\xi_i}{\eta_i}$ and $\fsnorm{f_i}=\norm{\xi_i}\norm{\eta_i}$, for some
representations $\theta_i:\Lm_i \to \cU(H_i)$ and vectors $\xi_i,\eta_i\in H_i$.
Put $H=H_1\tp_2 H_2$. Then $\theta_1\tp\theta_2 : \Lm_1\times\Lm_2\to \cU(H)$ is a representation and
\[ (f_1\tp f_2)(x,y) = \ip{(\theta_1\tp\theta_2)(x,y)(\xi_1\tp\xi_2)}{\eta_1\tp\eta_2}\;. \]
Hence $f_1\tp f_2\in \FS(\Lm_1\times\Lm_2)$ and
\[
\fsnorm{f_1\tp f_2} \leq \norm{\xi_1\tp\xi_2}_H \norm{\eta_1\tp\eta_2}_H \leq \norm{\xi_1} \norm{\xi_2} \norm{\eta_1}\norm{\eta_2} = \fsnorm{f_1} \fsnorm{f_2}\;.
\]
For the converse inequality, observe that for any
$b_1\in c_{00}(\Lm_1)$ and
$b_2\in c_{00}(\Lm_1)$,
\[ \begin{aligned}
\norm{b_1\tp b_2}_{\Cst(\Lm_1\times\Lm_2)} \fsnorm{f_1\tp f_2}
& \geq
\left\vert \sum_{(x,y)\in\Lm_1\times\Lm_2} (f_1\tp f_2)(x,y) (b_1\tp b_2)(x,y) \right\vert \\
&  = \left\vert \sum_{x\in\Lm_1} f_1(x) b_1(x) \sum_{y\in\Lm_2} f_2(y)b_2(y)\right\vert,
\end{aligned} \]
Therefore, if we can prove that
\begin{equation}\label{eq:cstar-tensor}
\norm{b}_{\Cst(\Lm_1\times\Lm_2)} \leq \norm{b_1}_{\Cst(\Lm_1)} \norm{b_2}_{\Cst(\Lm_2)}
\tag{$*$}
\end{equation}
it will follow that $\fsnorm{f_1\tp f_2}\geq \fsnorm{f_1}\fsnorm{f_2}$ as required.

Let $\theta:\Lm_1\times\Lm_2 \to \cU(H)$ be a representation and let $\theta_1(x)=\theta(x,e_2)$, $\theta_2(y)=\theta(e_1,y)$, where $e_1$ and $e_2$ denote the identity elements of $\Lm_1$ and $\Lm_2$ respectively. Then
\[\theta(b)
= \sum_{(x,y)\in \Lm_1\times\Lm_2} b_1(x)b_2(y) \theta(x,y)
= \sum_{(x,y)\in \Lm_1\times\Lm_2} b_1(x)b_2(y) \theta_1(x)\theta_2(y) = \theta_1(b_1)\theta_2(b_2)
\]
so that $\norm{\theta(b)} \leq \norm{b_1}_{\Cst(\Lm_1)} \norm{b_2}_{\Cst(\Lm_2)}$. Taking the supremum over all $\theta$ yields \eqref{eq:cstar-tensor}.
\end{proof}

The following result is surely not new, but since we did not find an explicit reference in the literature, we have included the statement and proof for sake of completeness.

\begin{lem}\label{l:FS countable}
Let $\Lm$ be a discrete group and let $f\in\FS(\Lm)$. Then there exists a countable subgroup $\Lm_0\leq\Lm$ such that $\fsnorm{{f\vert}_{\Lm_0}}=\fsnorm{f}$.
\end{lem}

\begin{proof}
We start with a more general observation: whenever $\Delta\leq\Lm$ we have $\fsnorm{f1_{\Delta}} = \fsnorm{{ f\vert}_\Delta}$. This is a special case of results for extending coefficient functions defined on open subgroups, see e.g.\ \cite[Lemma 7.1.2]{KL_AGBGbook} for a sketch of the proof.

Now choose a sequence $(b_n)$ in $c_{00}(\Lm)$ such that $\norm{b_n}_{\Cst(\Lm)} \leq 1$ and
\begin{equation}\label{eq:sandwich}
\fsnorm{f}- \frac{1}{n}
\leq \Abs{ \sum_{x \in \Lm} f(x)b_n(x) }  \leq \fsnorm{f}.
\tag{$*$}
\end{equation}
For each $n$ let $X_n\subset\Lm$ be a finite subset such that $\operatorname{supp}(b_n) \subseteq X_n$,
and let $\Lm_0$ be the subgroup of $\Lm$ generated by $\bigcup_{n\geq 1} X_n$, which is clearly countable.
The middle term in \eqref{eq:sandwich} is equal to $\Abs{ \sum_{x\in\Lm} (f1_{\Lm_0})(x) b_n(x) }$, which by definition of the norm in $\FS(\Lm)$ is bounded above by $\fsnorm{f1_{\Lm_0}}$. Hence $\fsnorm{f}\leq \fsnorm{f1_{\Lm_0}}$.
On the other hand, by our initial remarks and Lemma \ref{l:FS basics}\ref{li:restriction},
\[
\fsnorm{f1_{\Lm_0}} = \fsnorm{{f\vert}_{\Lm_0}} \leq \fsnorm{f}\;,
\]
and the desired result follows.
\end{proof}

We now recall some definitions and notation that were already used in the introduction.
Given a (locally compact) group $G$, the \dt{anti-diagonal of $G\times G$} is the set $\adiag(G)=\{ (x,x^{-1})\colon x\in G\}$.

\begin{dfn}
Given a discrete group $\Gm$ we define the \dt{anti-diagonal constant of $\Gm$} to be
\[
\AD(\Gm) = \norm{1_{\adiag(\Gm)}}_{\FS(\twice{\Gm})}
\]
with the convention that $\AD(\Gm)=\infty$ if $1_{\adiag(\Gm)} \notin \FS(\twice{\Gm})$.
\end{dfn}

As remarked in the introduction, Forrest and Runde showed that $\AD(\Gm)<\infty$ if and only if $\Gm$ is~VA. We will use this fact freely without further comment.

\begin{prop}\label{p:AD gen prop}
If $\Gm$ is a discrete group and $\Gm_0\leq\Gm$ then $\AD(\Gm_0)\leq\AD(\Gm)$.
If $\Gm_1$ and $\Gm_2$ are discrete groups then
$\AD(\Gm_1\times \Gm_2) =\AD(\Gm_1)\AD(\Gm_2)$.
\end{prop}

\begin{proof}
Both statements follow immediately from the corresponding parts of Lemma~\ref{l:FS basics}.
\end{proof}

We sometimes refer to the first statement in Proposition \ref{p:AD gen prop} as \dt{monotonicity} of~$\AD$; it will be our main tool for most of the paper, since it allows us to get lower bounds on $\AD(\Gm)$ by choosing an appropriate countable subgroup of $\Gm$. The next result is less essential, but has independent interest, and it will be convenient to make use of it in Section~\ref{s:classify minimal AD}.

\begin{prop}[Countable saturation for $\AD$]
\label{p:AD countable-sat}
Every VA group $\Gm$ has a countable subgroup $\Gm_0$ such that $\AD(\Gm)=\AD(\Gm_0)$.
\end{prop}

\begin{proof}
By Lemma \ref{l:FS countable} applied to $\FS(\twice{\Gm})$, there is a countable $\Lm_0 \leq \twice{\Gm}$ such that
$\fsnorm{1_{\adiag(\Gm)}} = \fsnorm{ {1_{\adiag(\Gm)}\vert}_{\Lm_0}}$.
Let $\Gm_0$ be the subgroup of $\Gm$ generated by the subsets $\{ x\colon (x,y)\in\Lm_0\}$ and $\{ y \colon (x,y)\in \Lm_0\}$. Then $\Gm_0$ is countable and $\Lm_0 \leq \Gm_0\times\Gm_0$, so that (by two applications of
Lemma~\ref{l:FS basics}\ref{li:restriction})
\[
\fsnorm{1_{\adiag(\Gm)}}
\geq \fsnorm{1_{\adiag(\Gm_0)}}
\geq \fsnorm{ {1_{\adiag(\Gm)}\vert}_{\Lm_0} }. 
\]
Hence $\fsnorm{1_{\adiag(\Gm)}}
=
\fsnorm{1_{\adiag(\Gm_0)}}$, as required.
\end{proof}

\end{subsection}

\begin{subsection}{Irreducible representations and the unitary dual}
Throughout this paper, irreducibility of a (unitary) representation is understood in the topological sense, i.e.\ the only closed invariant subspaces are $\{0\}$ and the whole space. For convenience, we will abbreviate the phrase ``irreducible representation'' to ``irrep''.

Given a locally compact group $G$, we write $\Ghat$ for the \dt{unitary dual} of $G$: that is, the set of unitary equivalence classes of continuous irreps of~$G$.
To simplify formulas and statements of results, for most of this paper we shall follow a standard abuse of notation, and identify $\Ghat$ with a set of representative irreps of $G$: that is, we tacitly fix a function $z\mapsto \pi_z$ with the property that $[\pi_z]=z$.
With this convention, already seen several times in Section~\ref{s:introduction}, we refer to particular irreps as being elements of $\Ghat$.

\begin{rem}\label{r:VA implies Type I}
If $\Gm$ is a VA group, then in fact it has a \emph{normal} abelian subgroup $\Lm$ such that $|\Gm:\Lm|<\infty$. Therefore, as shown in \cite[Satz 5]{Thoma68}, the following properties hold:
\begin{itemize}
\item Every factor representation of $\Gm$ is Type I (in other words, $\Gm$ is a \dt{Type~I group});
\item Every irrep of $\Gm$ is finite-dimensional and has degree bounded above by $|\Gm:\Lm|$.
\end{itemize}
This permits us to freely make use of the ``toolkit'' that has been developed for harmonic analysis on Type~I locally compact groups. For instance: if $\Gm$ is a \emph{countable} VA group, then it is second-countable and Type~I, so the Borel $\sigma$-algebra on $\Gmhat$ is standard in the sense of measure theory. (See \cite{DixCstar} for details, in particular Proposition 4.6.1 and Section 13.3.)
\end{rem}

In both the proof and the applications of the main formula (Equation \eqref{eq:AD-explicit-formula}), we will freely use results and ideas from the theory of the non-commutative Fourier transform and Fourier inversion.
We give an abridged review of the necessary notation and results from the literature, specialized to the setting of countable VA groups.
 For an overview of more general results, see \cite[Section 7.5]{Folland_ed2}.

Let $\Gm$ be a countable VA group.
Given $f\in c_{00}(\Gm)$,
note that when $\pi$ is a finite-dimensional representation of $\Gm$ the scalar $\Tr(\pi(f)\pi(f)^*$ depends only on the unitary equivalence class of~$\pi$. Hence each such $f$ yields a well-defined function $\Gmhat\to\Cplx$, given by $[\pi]\mapsto\Tr(\pi(f)\pi(f)^*)$, and one can show this function is continuous with respect to the Fell topology.
Furthermore, there is a non-negative (Radon) measure $\nu$ on $\Gmhat$, the \dt{Plancherel measure}, such that
\begin{equation}\label{eq:plancherel-measure}
\norm{f}_2^2 = \int_{\Gmhat} \Tr( \pi(f)\pi(f)^* )\,d\nu(\pi)
\qquad\text{for every $f\in c_{00}(\Gm)$.}
\end{equation}
Taking $f=\delta_e$ we have $1 = \int_{\Gmhat} d_\pi\,d\nu(\pi)$.
 (Note that we are tacitly adopting the normalization where $\Gm$ carries counting measure.)
Moreover, by combining \eqref{eq:plancherel-measure} with a polarization argument, we deduce that
\begin{equation}\label{eq:plancherel-identity}
{\ip{f}{g}}_{\ell^2(\Gm)} = \int_{\Gmhat} \Tr( \pi(f)\pi(g)^* )\,d\nu(\pi)
\qquad\text{for every $f,g\in c_{00}(\Gm)$.}
\end{equation}

Note that in Equations \eqref{eq:plancherel-measure} and \eqref{eq:plancherel-identity} we are following the convention mentioned earlier, treating $\Gmhat$ as a set of irreps by means of some tacit selection map. This practice will be more convenient when working with the main formula for $\AD(\Gm)$. However, in order to obtain this formula (i.e.\ in order to prove Theorem \ref{t:main formula}) we need to address measurability issues for certain operator fields defined on $\Gmhat$ and $\hatofsquare{\Gm}$, and here it seems more appropriate to make explicit this selection of an irrep from each equivalence class.

We therefore give a pr\'ecis of the required facts concerning measurable fields of Hilbert spaces, operators and representations. For all relevant definitions concerning these objects, the reader may consult the overview in \cite[Section 7.4]{Folland_ed2}, or the presentation in \cite[Part II, Ch.~1--3]{DixVNA} which has full details. The following terminology, however, is non-standard and made purely for convenience of later reference.

\begin{dfn}\label{d:selection-of-Ghat}
Let $G$ be a locally compact group and equip $\Ghat$ with the $\sigma$-algebra generated by the Fell topology. A \dt{measurable selection for $\Ghat$} is a pair $(\cH_\bullet,\beta_\bullet)$, where:
\begin{itemize}
\item for each $z\in\Ghat$, $\beta_z$ is a (continuous) representation of $G$ acting on a Hilbert space $\cH_z$, such that $\beta_z$ belongs to the equivalence class $z$;
\item $z\mapsto\cH_z$ is a measurable field of Hilbert spaces;
\item for each $g\in G$, $z\mapsto \beta_z(g)$ is a measurable operator field with respect to $\cH$.
\end{itemize}
\end{dfn}

It is not immediately obvious that we can always find measurable selections in the sense of the previous definition, even for countable VA groups.
The following result will be sufficient for our needs: it is the specialization to the countable VA setting of results, due originally to Mackey, that hold for second-countable Type I groups.

\begin{lem}[c.f.\ {\cite[Lemma 7.31]{Folland_ed2}}]\label{l:measurable selection}
Let $\Delta$ be a countable VA group, and equip $\Delhat$ with the $\sigma$-algebra generated by the Fell topology. For each $z\in\Delhat$, let $d(z)$ be the degree of any irrep belonging to the equivalence class $z$, and equip $\Cplx^{d(z)}$ with its standard inner product. Then $(\Cplx^{d(z)})_{z\in\Delhat}$ is a measurable field of Hilbert spaces on $\Delhat$, and there is a measurable field of representations on $\Delhat$, denoted by $z\mapsto \alpha_z$,
such that $\alpha_z$ acts on $\Cplx^{d(z)}$ and $\alpha_z \in z$ for each $z\in \Delhat$.
In particular, we have a measurable selection for~$\Delhat$.
\end{lem}

(The field $\Cplx^{d(\bullet)}$ is often referred to as the \dt{canonical field of Hilbert spaces} over~$\Delhat$.)

\end{subsection}

\begin{subsection}{The dual of the product of two countable VA groups}
\label{ss:dual-of-product}
To prove Theorem \ref{t:main formula}, when considering a countable VA group~$\Gm$ we need to work with measures and operator fields on (the Borel $\sigma$-algebra of) $\hatofsquare{\Gm}$, rather than on $\Gmhat$ itself. We also need to consider measurable selections for $\hatofsquare{\Gm}$ which are adapted to the product structure of $\twice{\Gm}$, and this is not automatic if Lemma \ref{l:measurable selection} is applied with $\Delta=\twice{\Gm}$.

Therefore, in this subsection we give an overview of known results which, in the specific setting of countable VA groups, allow us to identify $\hatofsquare{\Gm}$ with $\twice{\Gmhat}$ not just as sets but as measure spaces (a precise statement is given in Proposition \ref{p:identification} below).
Compared with the previous subsections, our presentation is more detailed, because we could not find explicit statements and proofs in \cite{DixCstar}, \cite{KL_AGBGbook} or \cite{KT_book}; our aim is to provide a sufficient road-map for those who wish to check details in full.

As motivation, recall how this works for finite groups. If $G_1$ and $G_2$ are finite groups, and $\pi_1$ and $\pi_2$ are irreps of $G_1$ and $G_2$ respectively, then $\pi_1\tp\pi_2$ is an irrep of $G_1\times G_2$, and we obtain a well-defined function
\begin{equation}\label{eq:define-J}
J:\hatonehattwo{G}  \to \hatonetwo{G}
\quad,\quad
([\pi_1],[\pi_2]) \mapsto [\pi_1 \otimes\pi_2]
\end{equation}
Since irreps of finite groups are classified by their traces, one can use the theory of group characters to show that $J$ is a bijection.

For general locally compact groups $G_1$ and $G_2$, the map $J$ from \eqref{eq:define-J} is still well-defined and injective.
A readable account is given\footnotemark\ in \cite[Section 7.3]{Folland_ed2}, with a minor caveat:
to prove injectivity of $J$, one uses the fact that if $\pi_1,\sigma_1$ are irreps of $G_1$ and $V$, $W$ are Hilbert spaces such that $\pi_1\tp I_V \sim \pi_2\tp I_W$, then $\pi_1\sim \pi_2$. This fact --- that irreps which are quasi-equivalent are equivalent --- does not follow directly from \cite[Proposition 7.13]{Folland_ed2} as claimed, but it can be found as \cite[Proposition 5.3.3]{DixCstar}. (It is also not hard to prove directly using Schur's lemma and basic properties of commutants.)

Moreover, $J:\hatonehattwo{G}\to \hatonetwo{G}$ is continuous with respect to the natural topologies on domain and codomain, which we might informally refer to as ``product-of-Fell'' and ``Fell-of-product''. This seems to be well known: see \cite[Theorem~2]{Fell63_wc-tp} or \cite[Proposition 5.3]{KT_book} for the proof of a more general result
(and see also Remark \ref{r:proposed by referee} below).

\begin{rem}\label{r:bijection for Type I}
In general $J$ may fail to be surjective. However, it is surjective (and therefore bijective) if either $G_1$ or $G_2$ is Type~I (see \cite[Theorem 7.17]{Folland_ed2} for a proof).
\end{rem}

We can now state the desired ``identification'' of measure spaces.

\begin{prop}\label{p:identification}
Let $\Gm_1$ and $\Gm_2$ be countable VA groups (which implies that $\Gm_1\times\Gm_2$ is also countable~VA). Let $\Sigma_1$, $\Sigma_2$ and $\Sigma_{12}$ be the $\sigma$-algebras on $\widehat{\Gm_1}$, $\widehat{\Gm_2}$ and $\hatonetwo{\Gm}$ generated by the Fell topologies on these spaces, and let $\Sigma_1\boxtimes\Sigma_2$ denote the product $\sigma$-algebra on $\hatonehattwo{\Gm}$.

\begin{romnum}
\item\label{li:standard}
The following measurable spaces are all standard (i.e.\ their $\sigma$-algebras are generated by some choice of Polish topology on the underlying space):
\[
(\widehat{\Gm_1},\Sigma_1)
\;;\;
(\widehat{\Gm_2},\Sigma_2)
\;;\;
(\hatonetwo{\Gm},\Sigma_{12})
\;;\;
(\hatonehattwo{\Gm},\Sigma_1\boxtimes\Sigma_2).
\]
\item\label{li:punchline}
 $J$ is a Borel isomorphism from $(\hatonehattwo{\Gm},\Sigma_1\boxtimes\Sigma_2)$ onto $(\hatonetwo{\Gm},\Sigma_{12})$.
\end{romnum}
\end{prop}

\begin{proof}
Each of the first three cases in part~\ref{li:standard} follows from the results mentioned in Remark~\ref{r:VA implies Type I}.
The fourth case follows from the fact that the product of standard measurable spaces is also standard. So only part \ref{li:punchline} requires further justification.

Denote the Fell topologies on $\widehat{\Gm_1}$, $\widehat{\Gm_2}$ and $\widehat{\Gm_{12}}$ by $\tau_1$, $\tau_2$ and $\tau_{12}$. Denote the product topology on $\hatonehattwo{\Gm}$ by~$\tau_1\boxtimes\tau_2$, and let $\Sigma(\tau_1\boxtimes\tau_2)$ denote the $\sigma$-algebra it generates.
By the remarks before this proposition, $J$ is $(\tau_1\times\tau_2)$-to-$\tau_{12}$-continuous, and so it is $\Sigma(\tau_1\boxtimes\tau_2)$-to-$\Sigma_{12}$-measurable.

Clearly $\Sigma(\tau_1\boxtimes\tau_2)\supseteq \Sigma_1\boxtimes\Sigma_2$ since the product $\sigma$-algebra is generated by ``open rectangles''. We now show that the converse inclusion holds. First, note that $\tau_1$ is second-countable: this is not mentioned explicitly in \cite{DixCstar}, but follows from countability of $\Gm_1$ and the fact that $(\widehat{\Gm_1},\tau_1)$ is a topological quotient of the pure state space of $\Cst(\Gm_1)$ with the relative \wstar-topology (see \cite[paragraph 3.4.12]{DixCstar}).
The same reasoning shows that $\tau_2$ is second-countable. Thus each $V\in \tau_1\boxtimes\tau_2$ is a \emph{countable} union of ``open rectangles'', so belongs to $\Sigma_1\boxtimes\Sigma_2$ as required.

From this, we see that $J$ is $\Sigma_1\boxtimes\Sigma_2$-to-$\Sigma_{12}$ measurable, where both $\sigma$-algebras are standard by part~\ref{li:standard}.
Also, by Remarks \ref{r:VA implies Type I} and \ref{r:bijection for Type I}, $J:\hatonehattwo{\Gm}\to \hatonetwo{\Gm}$ is bijective.  Therefore its inverse is also measurable by a theorem of Souslin (\cite[B22]{DixCstar}; for details see e.g.~\cite[Corollary~A.10]{Takesaki_vol1}).
\end{proof}

\begin{rem}\label{r:proposed by referee}
Our use of Souslin's theorem was made for convenience rather than necessity. In response to an earlier version of this paper, it was pointed out by the referee that the following result holds, without any assumptions of second-countability or the Type~I property:
\begin{quote}
\textit{for any locally compact groups $G_1$ and $G_2$, the canonical map $J:\hatonehattwo{G}\to \hatonetwo{G}$ is both continuous and relatively open.}
\end{quote}
 (Here, a \dt{relatively open mapping} between topological spaces is one which maps open subsets of the domain onto relatively open subsets of the image.)
In particular, in the proof of Proposition \ref{p:identification}, $J^{-1} : \hatonetwo{\Gm} \to \hatonehattwo{\Gm}$ is continuous, so measurability follows without requiring the appeal to Souslin's theorem.
\end{rem}

To the author's knowledge, the general result quoted in Remark \ref{r:proposed by referee} is not discussed in \cite{DixCstar}, \cite{Folland_ed2}, \cite{KL_AGBGbook} or \cite{KT_book}. Although it is not needed for the present paper, it might be useful for future work, so we have provided an expanded version of the referee's sketch in the appendix.

\paragraph{Author's note.} After this paper was accepted for publication, I learned that the main results of Appendix \ref{app:general-dual-of-product} can also be found in \cite[Appendix B.5]{RW_morita-book}, with similar proofs. Since this is not a particularly obvious place to look for such results, and since it is remarked at the start of \cite[Appendix B.5]{RW_morita-book} that ``it is not easy to dig detailed proofs of the resulting theorems out of the literature'', I have decided to keep Appendix~\ref{app:general-dual-of-product} for the reader's convenience.
\end{subsection}

\end{section}


\begin{section}{Fourier coefficients of the anti-diagonal}

In this section we will prove Theorems \ref{t:AD of finite}, \ref{t:main formula} and \ref{t:sharp LB}.

For a compact group $K$ and $h\in\FA(K)$, the Fourier norm of $h$ can be calculated from its matrix-valued Fourier coefficients.
We only need the case when $K$ is a finite group; with our earlier convention that $\sigma(h)=\sum_{x\in K} h(x)\sigma(x)$, we obtain
\begin{equation}\label{eq:fourier-norm-compact}
\norm{h}_{\FA(K)} = \sum_{\sigma\in \widehat{K}} \frac{d_\sigma}{|K|} \norm{\sigma(h)}_{(1)} \;.
\end{equation}
In particular, for finite $G$ we can calculate $\AD(G)$ by working out the Fourier coefficients of $1_{\adiag(G)}$, which turn out to have an explicit and tractable form if one exploits Schur orthogonality relations for $G$ itself.

The calculation of the Fourier coefficients works just as well as for an infinite compact~$G$, provided that one regards $1_{\adiag(G)}$ as a singular measure supported on the anti-diagonal.
%
Given a Radon measure $\gm$ on a compact group $K$, we define its matrix-valued Fourier coefficients $(\pi(\gm))_{\pi\in\widehat{K}}$ as integrals in the weak sense, i.e.
\[ \ip{\pi(\gm)\xi}{\eta} \defeq \int_K \ip{\pi(x)\xi}{\eta}\,d\gm(x) \qquad(\xi,\eta\in H_\pi) \]
For $K$ finite, we may view $\gm$ as a function on $K$, and then this definition of $\pi(\gm)$ agrees with the earlier notation used in Section~\ref{ss:AG-BG} for finitely supported functions on discrete groups.

Given a compact group $G$ equipped with a choice of Haar measure $\mu$, we introduce a Radon measure $\muad$ on $G\times G$,
which is defined to be the image measure of~$\mu$ (also called the pushforward measure of $\mu$) with respect to the embedding $\adiag:G\to G\times G$. Explicitly,
\begin{equation}\label{eq:dirac on AD}
\int_{G\times G} f(x,y)\, d\muad(x,y) \defeq \int_G f(x,x^{-1})\,d\mu(x) \qquad(f\in C(G\times G)).
\end{equation}
Note that if $G$ is finite and $\mu$ is counting measure, then the Radon--Nikodym derivative of $\muad$ with respect to counting measure on $G\times G$ is~$1_{\adiag(G)}$.

\begin{rem}\label{r:finite case}
If $G$ is finite, but we chose to equip $G$ and $G\times G$ with their uniform measures, as is more common when discussing Fourier analysis on compact groups, then the Radon-Nikodym derivative of $\muad$ with respect to uniform measure on $G\times G$ would be $|G| 1_{\adiag(G)}$.
\end{rem}

\begin{dfn}
Let $H$ be a Hilbert space. The \dt{flip map} on $H\tp_2 H$ is defined to be the unique bounded operator $\sX_H : H\tp_2 H \to H\tp_2 H$ which satisfies $\sX_H(\xi\tp\xi') = \xi'\tp\xi$ for every $\xi,\xi'\in H$.
\end{dfn}

\begin{prop}[Fourier coefficients of the anti-diagonal]\label{p:compact FC of AD}
Let $G$ be compact and let $\mu$, $\muad$ be as above.
Let $\pi,\sigma\in\Ghat$. Then
\begin{equation}\label{eq:FC-of-antidiag}
(\pi\tp\sigma)(\muad) = \begin{cases} 0 & \text{if $\pi\not\sim\sigma$,} \\ d_\pi^{-1}\mu(G)\sX_{\pi} & \text{if $\pi=\sigma$,} 
\end{cases}
\end{equation}
where by mild abuse of notation $\sX_\pi$ denotes the flip map on $H_\pi$.
\end{prop}

\begin{proof}
By linearity and continuity it suffices to check that the identities in \eqref{eq:FC-of-antidiag} hold when both sides are evaluated on elementary tensors in $H_\pi\tp H_\sigma$.
Let $\xi,\eta\in H_\pi$ and $\xi',\eta'\in H_\sigma$. Then
\[ \begin{aligned}
\ip{ (\pi\tp\sigma)(\muad)(\xi\tp\xi')}{\eta\tp\eta'}
& =
\int_{G\times G} \ip{\pi(x)\xi}{\eta} \ip{\sigma(y)\xi'}{\eta'} \,d\muad(x,y) \\
& =
\int_{G} \ip{\pi(x)\xi}{\eta} \ip{\sigma(x^{-1})\xi'}{\eta'} \,d\mu(x) \\
& =
\int_{G} \ip{\pi(x)\xi}{\eta} \overline{\ip{\sigma(x)\eta'}{\xi'}} \,d\mu(x) \\
\end{aligned} \]
By the Schur orthogonality relations for $G$, the right-hand side vanishes if $\pi\not\sim\sigma$, while for $\pi=\sigma$ we have
\[ \begin{aligned}
\int_{G} \ip{\pi(x)\xi}{\eta} \overline{\ip{\pi(x)\eta'}{\xi'}} \,d\mu(x)
& = \frac{\mu(G)}{d_\pi} \ip{\xi'}{\eta} \ip{\xi}{\eta'}
& = \frac{\mu(G)}{d_\pi} \ip{\sX_\pi(\xi\tp\xi')}{\eta\tp\eta'}
\end{aligned} \]
as claimed.
\end{proof}

\begin{rem}\label{r:fin-dim-flip}
When $\dim(H)<\infty$ there is an explicit formula for $\sX_H$ with respect to a choice of o.n.\ basis for $H$ and the ensuing o.n.\ basis for $H\tp_2 H$ (in fact, this was how the formulas in Proposition~\ref{p:compact FC of AD} were originally obtained). 
Let $(E_{ij})$ be the standard matrix units for $M_n(\Cplx)$ and consider $\sum_{i,j=1}^n E_{ij} \tp E_{ji}$ viewed as an operator on $\Cplx^n \tp \Cplx^n$. A direct calculation shows that this operator implements the flip map, and hence it equals $\sX_{\Cplx^n}$.

The referee has pointed out that this formula can be extended to an arbitrary Hilbert space~$H$: fixing a maximal orthonormal family $(e_i)_{i\in {\mathbb I}}$ in $H$ and defining corresponding ``matrix units'' $E_{ij}\in\Bdd(H)$, we have $\sX_H = \sum_{i,j\in {\mathbb I}} E_{ij} \tp E_{ji}$ where the infinite sum denotes unordered convergence (over finite subsets of ${\mathbb I}$) in the strong operator topology of $\Bdd(H\tp_2 H)$.
\end{rem}

\begin{proof}[The proof of Theorem~\ref{t:AD of finite}]
It suffices to prove the first equality in \eqref{eq:AD of finite}, since the second one is Johnson's formula \eqref{eq:BEJ formula}. Let $G$ be a finite group. Then $\FS(G\times G) = \FA(G\times G)$, so applying Equation~\eqref{eq:fourier-norm-compact} with $K=\twice{G}$ yields
\[
\AD(G)
 = \norm{1_{\adiag(G)}}_{\FA(G\times G)} 
 = \sum_{(\pi,\sigma) \in \twice{\Ghat}} \frac{d_{\pi\tp\sigma}}{|G\times G|} \norm{ (\pi\tp\rho)(1_{\adiag(G)})}_{(1)} \;.
 \]
Note that we may identify $\hatofsquare{G}$ with $\twice{\Ghat}$; and if $\pi,\sigma\in\Ghat$, Proposition \ref{p:compact FC of AD} gives
\[
(\pi\tp\sigma)(1_{\adiag}(G)) = \begin{cases} 0 & \text{if $\pi\not\sim\sigma$,} \\ d_\pi^{-1}|G|\sX_{\pi} & \text{if $\pi=\sigma$,} 
\end{cases}
\]
and so
\[
\AD(G) = \sum_{\pi\in\Ghat}\frac{(d_\pi)^2}{|G|^2} \cdot\frac{|G|}{d_\pi} \norm{\sX_\pi}_{(1)} \;.
\]
Note that $\sX_\pi$ is a Hermitian involution on $H_\pi\tp_2 H_\pi$, and so with respect to a suitable o.n.\ basis of $H_\pi\tp_2 H_\pi$ it is represented by a diagonal matrix with $\pm 1$ entries. Hence $\norm{\sX_\pi}_{(1)} = \dim(H_\pi \tp_2 H_\pi) = (d_\pi)^2$, and the desired identity for $\AD(G)$ follows.
\end{proof}

While Proposition~\ref{p:compact FC of AD} does not apply to infinite discrete groups, it does suggest how we might proceed for groups which have a well-behaved notion of Fourier transform and Fourier inversion.
However, since $1_{\adiag(\Gm)}$ does not belong to $\FA(\twice{\Gm})$, we cannot calculate its Fourier--Stieltjes norm by a direct application of the Plancherel theorem for $\twice{\Gm}$. Instead, our strategy is to exploit results of Arsac, developed in \cite[Section~3G]{Ars_Lyon76}, which provide an isometric identification between certain subspaces of $\FS(\twice{\Gm})$ and spaces of ``operator-valued measures'' on $\hatofsquare{\Gm}$.
%
The following is a special case of \cite[(3.55)]{Ars_Lyon76}.
\begin{thm}[Arsac]\label{t:arsac-bochner}
Let $G$ be a second-countable Type I group, and fix a measurable selection $(\cH_\bullet,\beta_\bullet)$ for $\Ghat$.
If $\mu$ is a $\sigma$-finite Borel measure on $\Ghat$, and $T=(T_z)_{z\in\Ghat}$ is a $\mu$-measurable operator field acting on $\cH$, such that
\[
\norm{T}_{1,\mu} \defeq \int_{\Ghat} \norm{T_z}_{(1)}\,d\mu(z) <\infty,
\]
then the function $\Psi_\mu(T):G\to\Cplx$ defined by
\[
\Psi_\mu(T)(s) \defeq \int_{\Ghat} \Tr(T_z\beta_z(s))\,d\mu(z)
\]
belongs to $\FS(G)$ and satisfies $\fsnorm{\Psi_\mu(T)}= \norm{T}_{1,\mu}$.
\end{thm}

\begin{rem}\label{r:extra-structure}
The drawback of appealing to Theorem \ref{t:arsac-bochner} as a black box, is that both the proof and statement of that theorem ignore the extra structure present when $G$ has the form $\twice{\Gm}$. The arguments in \cite[Section 3G]{Ars_Lyon76} use the stratification of $\Ghat$ according to the degrees of irreps, but for a product of two groups $G=G_1\times G_2$ this stratification does not ``see'' how an irrep of $G$ might decompose as a tensor product of irreps of $G_1$ and $G_2$, and in the case where $G_1=G_2$ the stratification does not detect the ``diagonal subset'' inside $\Ghat$. (This is why we require something like Proposition \ref{p:identification} in our approach.)
\end{rem}

We will also need the following property of the flip map, which we state as a separate lemma since it has independent interest.
It can be derived  using the formula in Remark \ref{r:fin-dim-flip}, but we shall give a co-ordinate free proof.

\begin{lem}\label{l:trace-flip}
Let $B$ and $C$ be trace-class operators on a Hilbert space $H$. Then
\[
\Tr(\sX_H \cdot (B\tp C) ) = \Tr(BC).
\]
\end{lem}

\begin{proof}
By linearity and continuity it suffices to prove this when $B$ and $C$ are rank-one operators on $H$, say $B(\xi) = \ip{\xi}{\beta_2}\beta_1$ and $C(\eta)=\ip{\eta}{\gm_2}\gm_1$ for some $\beta_1,\beta_2,\gm_1,\gm_2\in H$. Then $BC(\eta) = \ip{\eta}{\gm_2}\ip{\gm_1}{\beta_2}\beta_1$, and so 
\[
\Tr(\sX_H\cdot (B\tp C))
= \ip{ \sX_H(\beta_1\tp\gm_1)}{\beta_2\tp\gm_2}
= \ip{\gm_1}{\beta_2} \ip{\beta_1}{\gm_2} = \Tr(BC)
\]
as required.
\end{proof}

\begin{proof}[The proof of the main formula (Theorem~\ref{t:main formula})]
The idea of the proof is to construct an explicit operator field $T$ over $\hatofsquare{\Gm}$, and a Borel measure $\nud$ on $\hatofsquare{\Gm}$, such that in the notation of Theorem \ref{t:arsac-bochner} with $G=\twice{\Gm}$ we obtain $\Psi_{\nud}(T)=1_{\adiag(\Gm)}$.

Let $(\Cplx^{d(\bullet)},\alpha_\bullet)$ be the measurable selection for $\Gmhat$ that is provided by Lemma \ref{l:measurable selection}.
By Proposition \ref{p:identification} we may identify $\hatofsquare{\Gm}$, equipped with the Borel $\sigma$-algebra generated by the Fell topology, with the product measure space $\twice{\Gmhat}$ (we suppress explicit mention of the bijection $J$ in order to simplify notation). With this identification made, routine calculations show that $(\Cplx^{d(\bullet)}\tp_2\Cplx^{d(\bullet)}, \alpha_\bullet\tp\alpha_\bullet)$ is a measurable selection for $\hatofsquare{\Gm}$.

Moreover, since $\Gmhat$ is a standard measure space, the diagonal subset of $\twice{\Gmhat}$ is measurable, and the diagonal embedding $\diag:\Gmhat\to\twice{\Gmhat}$ is a measurable function. Hence, the following constructions are well-defined:
\begin{itemize}
\item the image measure on $\twice{\Gmhat}$ obtained by pushing forward Plancherel measure $\nu$ on $\Gmhat$ along $\diag$ --- this is a Borel measure $\nud$ on $\twice{\Gmhat}$, satisfying
\[
\int_{\twice{\Gmhat}} h(x,y)\,d\nud(x,y) = \int_{\Gmhat} h(y,y)\,d\nu(y)
\]
for every bounded measurable function~$h$;
\item
a measurable operator field $T=(T_{x,y})_{(x,y)\in\twice{\Gmhat}}$ acting on $\Cplx^{d(\bullet)}\tp_2 \Cplx^{d(\bullet)}$, defined by
$T_{y,y} = \sX_{d(y)}$ and $T_{x,y}=0$ if $x\neq y$.
\end{itemize}
Here $\sX_d$ denotes the flip map on $\Cplx^d\tp_2\Cplx^d$. As seen earlier, $\norm{\sX_d}_{(1)}=d^2$.
Since $\nud$ is a finite measure and $\maxdeg(\Gm)<\infty$, this shows that $T$ is integrable in the sense required by Theorem \ref{t:arsac-bochner} (we shall give a more precise calculation later).

Now, given $s,t\in\Gm$, and using our identification of $\hatofsquare{\Gm}$ with $\twice{\Gmhat}$,
\[
\begin{aligned}
\Psi_{\nud}(T)(s,t)
 & = \int_{\twice{\Gmhat}} \Tr(T_{x,y} (\alpha_x\tp\alpha_y)(s,t)\,d\nud(x,y)
   & \text{(definition of $\Psi(T)$)} \\
 & = \int_{\Gmhat} \Tr\left(\sX_{d(y)} \cdot (\alpha_x(s)\tp\alpha_y(t))\right)\,d\nu(y)
   & \text{(definitions of $\nud$ and $T$)} \\
 & = \int_{\Gmhat} \Tr(\alpha_y(s)\alpha_y(t))\,d\nu(y)
   & \text{(Lemma \ref{l:trace-flip})} \\
 & = \ip{\delta_s}{\delta_{(t^{-1})}}
   & \text{(Equation \eqref{eq:plancherel-identity})}. \\
\end{aligned}
\]
Thus $\Psi_{\nud}(T) = 1_{\adiag(\Gm)}$. Hence, by Theorem \ref{t:arsac-bochner}, $\AD(\Gm)=\norm{T}_{1,\nud}$, which in turn is equal to
\[
\int_{\Delhat} \norm{T_z}_{(1)} \,d\nud(z) 
 = \int_{\Gmhat} \norm{\sX_{d(y)}}_{(1)} \,d\nu(y)
 = \int_{\Gmhat} d(y)^2 \,d\nu(y)
\]
as required.
\end{proof}

\begin{rem}\label{r:apologia}
It seems plausible that Theorem \ref{t:main formula} could be extended to cover uncountable VA groups, but to the author's knowledge, the only way to show this would be to develop non-separable analogues of the results in \cite[3G]{Ars_Lyon76}. This is rather annoying since Arsac's results, as stated, rely on the classical disintegration results listed in \cite{DixCstar}, which in turn are built on results for fields of von Neumann algebras that are usually developed under countability assumptions. It is probably also necessary to develop a version of Theorem~\ref{t:arsac-bochner} which is adapted to groups with a product decomposition, to avoid dealing with the countability assumptions  used in Section \ref{ss:dual-of-product}.
While not necessarily difficult, in total this work would be laborious, with limited rewards since Theorem~\ref{t:main formula} is most useful for countable groups anyway.
We have therefore chosen to avoid such technicalities in this paper.
\end{rem}

Let us now see how the main formula can be applied to prove the sharp lower bound (Theorem \ref{t:sharp LB}).

\begin{prop}\label{p:at-most-half}
Let $\Gm$ be a countable VANA group and let $\Om_1=\{ \pi \in\Gmhat \colon d_\pi =1\}$. Let $\nu$ be Plancherel measure on $\Gmhat$, normalized so that $\int_{\Gmhat}d_\pi\,d\nu(\pi)=1$. Then $\nu(\Om_1)\leq 1/2$.
\end{prop}

\begin{proof}
Pick $x$ and $y$ in $\Gm$ with $xy\neq yx$, and note that whenever $\pi\in\Om_1$, $\pi(xy)=\pi(yx)$ is a complex number of modulus~$1$.
The Plancherel formula
(Equation \ref{eq:plancherel-identity}) then gives us
\[ 0 = {\ip{\delta_{xy}}{\delta_{yx}}}_{\ell^2(\Gm)} = \int_{\Gmhat} \Tr( \pi(xy)\pi(yx)^* ) \,d\nu(\pi)
= \nu(\Om_1) + \int_{\Gmhat\setminus\Om_1} \Tr(\pi(xy)\pi(yx)^*))\,d\nu(\pi)
\]
and so
\[ \begin{aligned}
\nu(\Om_1)
 \leq \Abs{\int_{\Gmhat\setminus\Om_1} \Tr(\pi(xy)\pi(yx)^*))\,d\nu(\pi) }
& \leq \int_{\Gmhat\setminus\Om_1} d_\pi \, d\nu(\pi)   \\
& = \int_{\Gmhat} d_\pi\,d\nu(\pi) - \int_{\Om_1} d_\pi \,d\nu(\pi) \\
& = 1- \nu(\Om_1).
\end{aligned} \]
Rearranging gives $2\nu(\Om_1) \leq 1$ as required.
\end{proof}

\begin{proof}[The proof of Theorem \ref{t:sharp LB}]
The second statement follows from the first statement and the fact that $\AM(\FA(G))\geq \AD(G_d)$ for any locally compact group. Hence it suffices to prove that $\AD(\Gm)\geq 3/2$ whenever $\Gm$ is non-abelian.
	
Since $\AD(\Gm)=+\infty$ if $\Gm$ is not virtually abelian, we may assume without loss of generality that $\Gm$ is VANA.
Pick $x,y\in\Gm$ which do not commute, and let $\Gm_0$ be the subgroup of $\Gm$ generated by $x$ and $y$. Then $\Gm_0$ is countable and VANA.
 For each $n\in\Nat$, let $\Om_n = \{ \pi\in \widehat{\Gm_0} \colon d_\pi = n\}$. Then by Theorem~\ref{t:main formula}
\begin{equation}\label{eq:crude-lower-bound}
\AD(\Gm_0) - \nu(\Om_1) = \sum_{n\geq 2} n^2 \nu(\Om_n) \geq \sum_{n\geq 2} 2n\nu(\Om_n) = 2(1 - \nu(\Om_1).
\end{equation}
Rearranging and using Proposition \ref{p:at-most-half} yields $\AD(\Gm_0) \geq 2 - 1/2 = 3/2$. By monotonicity of $\AD$ (Proposition \ref{p:AD gen prop}) we deduce that $\AD(\Gm)\geq \AD(\Gm_0)\geq 3/2$, as required,
\end{proof}

In the next section we will characterize those non-abelian groups for which our lower bound on $\AD$ is sharp. To motivate some of the techniques, we finish this section by returning to an earlier example and calculating its anti-diagonal constant.

\begin{eg}\label{eg:AD red-heis-mod2}
Let $\Heis_2$ be the group from Example~\ref{eg:red-heis-mod2}. Recall that $\Heis_2$ is non-abelian yet has no finite non-abelian subgroups; in particular it is not of the form (abelian)$\times$(finite). Despite this, we shall see that it attains the lower bound in Theorem \ref{t:sharp LB}.

The irreps of $\Heis_2$ can be found by a standard application of the Mackey machine, with details in e.g.\ \cite[Section 6.8]{Folland_ed2}. We only require the following facts:
\begin{itemize}
\item every such representation is either $1$-dimensional or $2$-dimensional;
\item if $\pi$ is a $2$-dimensional irrep of $\Heis_2$, then~$\pi(0,0,[1])=-I_2$.
\end{itemize}
Let $\nu$ be Plancherel measure on $\widehat{\Heis_2}$.
The first property implies that $1= \nu(\Om_1) + 2 \nu(\Om_2)$.
If we put $x=(1,0,[0])$ and $y=(0,1,[0])$, direct calculation shows that $[x,y]=xy(yx)^{-1}=(0,0,[1])$.
Arguing as in the proof of Proposition~\ref{p:at-most-half}, we obtain
\[
0 = \ip{\delta_{xy}}{\delta_{yx}} = \int_{\pi \in\widehat{\Heis_2}} d_\pi \Tr(\pi([x,y])) \,d\nu(\pi)
= \nu(\Om_1) - 2\nu(\Om_2).
\]
Therefore $\nu(\Om_1)=1/2$ and $\nu(\Om_2)=1/4$, which yields $\AD(\Heis_2)=\nu(\Om_1)+4\nu(\Om_2) = 3/2$.
\end{eg}

\end{section}


\begin{section}{Characterizing groups with anti-diagonal constant $3/2$}
\label{s:classify minimal AD}

In this section we show that the groups $\Gm$ satisfying $\AD(\Gm)=3/2$ are exactly those satisfying $\Zindex{\Gm}=4$ (Theorem~\ref{t:minimizers of AD}). Note that we do not impose any countability restriction on~$\Gm$. However, since our ``main formula'' for the anti-diagonal constant is only established for countable VA groups, our approach inevitably passes through the countable setting.

Although the natural approach is to start with the condition $\AD(\Gm)=3/2$ and see what can be deduced, it is convenient to set up some general tools that can be used repeatedly.

\begin{dfn}\label{d:comm}
For a group $\Gm$ let $\comm(\Gm)=\{[x,y] \colon x,y\in\Gm\}$.
\end{dfn}

Note that $e\in\comm(\Gm)$, and $|\comm(\Gm)|=1$ if and only if $\Gm$ is abelian.
Thus the groups $\Gm$ for which $|\comm(\Gm)|=2$ are, informally, the non-abelian ones where there are as few commutators as possible. This has strong structural consequences; in particular, it follows from Lemma \ref{l:was-lem-C}\ref{li:was-C1} below that such groups are nilpotent of class~$2$.

\begin{lem}\label{l:was-lem-C}
Let $\Gm$ be a group satisfying $|\comm(\Gm)|=2$, and let $z$ be the non-identity element of $\comm(\Gm)$.
\begin{romnum}
\item\label{li:was-C1}
 $z^2=e$ and $z\in Z(\Gm)$.
\item\label{li:was-C2}
Let $\pi\in\Gmhat$ with $d_\pi>1$.  Then $\pi(z)=-I_\pi$. Hence, if $x$ and $y$ are noncommuting elements of $\Gm$, the operators $\pi(x)$ and $\pi(y)$ anti-commute, that is $\pi(x)\pi(y)+\pi(y)\pi(x)=0$.
\end{romnum}
\end{lem}

\begin{proof}\
If $\varphi$ is an automorphism of $\Gm$ then it maps $\comm(\Gm)\setminus\{e\}$ to itself. Thus $z$ is fixed by every automorphism of $\Gm$, in particular by all inner automorphisms, so $z\in Z(\Gm)$.
Since $([a,b])^{-1}=[b,a]$ for every $a,b\in \Gm$, we have $z^{-1}\in \comm(\Gm)\setminus\{e\}=\{z\}$. That is, $z^{-1}=z$. This proves part \ref{li:was-C1}.

Let $\pi\in\Gmhat$ with $d_\pi > 1$. Since $z\in Z(\Gm)$ and $\pi$ is irreducible, $\pi(z)\in\Cplx I_\pi$ by Schur's lemma. Moreover, $\pi(z)=\pm I_\pi$ since $z^2=e$. 
Suppose $\pi(z)=I_\pi$\,: then since $\comm(\Gm)=\{e,z\}$ we have $[\pi(a),\pi(b)] = \pi([a,b])=I_\pi$ for all $a,b\in\Gm$. But this implies that $\pi(\Gm)$ is abelian which contradicts irreducibility of~$\pi$.
Therefore $\pi(z)=-I_\pi$, and since $\pi(z)=\pi(x)\pi(y)\pi(x)^{-1}\pi(y)^{-1}$, rearranging gives $\pi(x)\pi(y)=-\pi(y)\pi(x)$. This proves part \ref{li:was-C2}.
\end{proof}

The relevance of the condition $|\comm(\Gm)|=2$ to the study of $\AD(\Gm)$ is not immediately apparent, but is partially explained by the next two results. In both results we restrict attention to countable VA groups, to make use of the Plancherel identity \eqref{eq:plancherel-identity}.

\begin{prop}\label{p:was-E}
Let $\Lm$ be a countable VA group with Plancherel measure $\nu$. Let $\Om_1=\{\pi\in\Lmhat \colon d_\pi=1\}$. Then $|\comm(\Lm)|=2$ if and only if $\nu(\Om_1)=1/2$.
\end{prop}

\begin{proof}
We may assume that $\Lm$ is non-abelian (since this follows from either of the assumptions $|\comm(\Lm)|=2$ or $\nu(\Om_1)=1/2$).

Note that $\Re\Tr(V)\leq d_\pi$ for any unitary $V\in\cU(H_\pi)$, with equality holding if and only if $V=I_\pi$. Hence, as already observed in the proof of Proposition \ref{p:at-most-half}:
\[
\begin{aligned}
\text{for all $z\in\comm(\Lm)\setminus\{e\}$,}
\quad \nu(\Om_1)
 & = \int_{\Lmhat\setminus\Om_1} -\Tr\pi(z)\,d\nu(\pi) \\
 & = \int_{\Lmhat\setminus\Om_1} \Re\Tr(-\pi(z))\,d\nu(\pi) \\
 &  \leq \int_{\Lmhat\setminus\Om_1} d_\pi \,d\nu(\pi) 
 & = 1- \nu(\Om_1),
\end{aligned}
\]
and equality holds if and only if $\pi(z)=-I_\pi$ for $\nu$-a.e.\ $\pi\in\Lmhat\setminus\Om_1$.

Note that $\nu(\Om_1)=1/2$ if and only if $\nu(\Om_1)=1-\nu(\Om_1)$. By the previous remarks, the following statements are therefore equivalent:
\begin{romnum}
\item $\nu(\Om_1)=1/2$;
\item\label{li:TWO}
$\pi(z)=-I_\pi$ for every $z\in \comm(\Lm)\setminus\{e\}$ and $\nu$-a.e.\ $\pi\in\Lmhat\setminus\Om_1$.
\end{romnum}

If $|\comm(\Lm)|=2$, then \ref{li:TWO} holds by  Lemma \ref{l:was-lem-C}\ref{li:was-C2}. On the other hand, suppose that \ref{li:TWO} holds.
Let $z,z'\in \comm{\Lm}\setminus\{e\}$. By assumption $\pi(z)=\pi(z')$ for $\nu$-a.e.~ $\pi\in\Lmhat\setminus\Om_1$; and for every $\pi\in\Om_1$ we have $\pi(z)=\pi(z')$. By injectivity of the non-abelian Fourier transform (or, alternatively, by applying the Plancherel identity \eqref{eq:plancherel-identity} to $\ip{\delta_z}{\delta_{z'}}$) it follows that $z=z'$. Thus $|\comm(\Lm)|=2$.
\end{proof}

\begin{lem}\label{l:was-A}
Let $\Lm$ be a countable VA group and let $\Om_n=\{\pi\in\Lmhat\colon d_\pi=n\}$. Let $\nu$ denote Plancherel measure on $\Lmhat$. The following are equivalent:
\begin{romnum}
\item\label{li:was-A1} $\AD(\Lm)=3/2$;
\item\label{li:was-A2} $\nu(\Om_1)=1/2$ and $\nu(\Om_n)=0$ for all $n\geq 3$.
\end{romnum}
\end{lem}

\begin{proof}
We already observed in the proof of Theorem \ref{t:sharp LB} --- specifically, in the inequality \eqref{eq:crude-lower-bound} --- that $\AD(\Lm)\geq 2-\nu(\Om_1)$, with equality if and only if $n^2\nu(\Om_n)=2n\nu(\Om_n)$ for all $n\geq 2$. In turn, this is equivalent to requiring $\nu(\Om_n)=0$ for all $n\geq 3$.

Now, Proposition \ref{p:at-most-half} implies that $2-\nu(\Om_1)\geq 3/2$ with equality if and only if $\nu(\Om_1)=1/2$. Thus \ref{li:was-A1} and \ref{li:was-A2} are equivalent.
\end{proof}

Let us sketch where we are headed. Lemma \ref{l:was-A} and Proposition \ref{p:was-E} indicate that, at least for countable groups, the assumption $\AD(\Lm)=3/2$ places strong restrictions on both the set of commutators in $\Lm$ and the set of possible degrees of irreps of $\Lm$. These restrictions turn out to imply that for any non-central $x\in\Lm$, the centralizer of $x$ has index $2$ in $\Lm$ but also contains $Z(\Lm)$ as a subgroup of index $2$; this will immediately imply that $\Zindex{\Lm}=4$. In the converse direction, if $\Gm$ is any group (not necessarily countable), assuming that $\Zindex{\Gm}=4$ turns out to have strong structural consequences for $\Gm$, which can be used to control both $\comm(\Gm)$ and the degrees of irreps of $\Gm$.

For technical reasons, let us switch attention from the assumption $\AD(\Gm)=3/2$ and its consequences to the assumption $\Zindex{\Gm}=4$ and its consequences. The following basic group-theoretic lemma can be found as a standard exercise in various texts.

\begin{lem}\label{l:was-prop-F}
Let $\Delta$ be a non-abelian group. Then $\Delta/Z(\Delta)$ is not cyclic. In particular $\Zindex{\Delta}\geq 4$.
\end{lem}

\begin{proof}[Outline of a proof]
The first part follows from a stronger statement: if $\Delta$ is a group and $N$ is a normal subgroup contained in $Z(\Delta)$ such that $\Delta/N$ is cyclic, then $\Delta$ is abelian. To prove this, lift the generator of $\Delta/N$ to some $r\in\Delta$, and note that $\Delta = \bigcup_{n\in\Zahl} r^n N \subseteq \bigcup_{n\in\Zahl} r^n Z(\Delta)$.
For the second part, it suffices to note that all groups of order $\leq 3$ are cyclic.
\end{proof}

The next result characterizes the groups in which the centre has index $4$, as those non-abelian groups where the set of commutators and the set of degrees of irreps are as small as possible. We have not found it stated in the literature, so we provide a detailed proof, since it may have interest independent of the application to anti-diagonal constants. Note that we do not require any countability assumptions.

\begin{thm}\label{t:was-G}
Let $\Gm$ be a group. The following conditions are equivalent.
\begin{romnum}
\item\label{li:was-G1}
$\Zindex{\Gm}=4$.
\item\label{li:was-G2}
$|\comm(\Gm)|=2$ and $\maxdeg(\Gm)=2$.
\item\label{li:was-G3}
$|\comm(\Gm)|=2$ and $\{\pi\in \Gmhat\colon d_\pi\leq 2\}$ separates the points of $\Gm$.
\end{romnum}
\end{thm}

\begin{rem}
The equivalence of parts \ref{li:was-G2} and \ref{li:was-G3} is a special case of  general results for the $\maxdeg$ function, which follow in turn from results on polynomial identities for noncommutative algebras. However, this does not seem to significantly shorten the proof of the equivalence \ref{li:was-G3}$\implies$\ref{li:was-G1}, which is what we need in the application to $\AD(\Gm)$.
\end{rem}

\begin{proof}[Proof of Theorem \ref{t:was-G}\ref{li:was-G1}$\implies$\ref{li:was-G2}]
Note that the commutator map $\twice{\Gm}\to\Gm$ factors through $(\Gm/Z(\Gm))\times (\Gm/Z(\Gm))$. More precisely, let $q:\Gm\to\Gm/Z(\Gm)$ be the quotient homomorphism, and define an equivalence relation $\sim$ on $\Gm$ by $x\sim y \iff q(x)=q(y)$. Given $x\sim x'$ and $y\sim y'$, a direct calculation shows that $[x,y]=[x',y']$.
Therefore, if $T$ is a transversal for $Z(\Gm)$ in $\Gm$, we have $\comm(\Gm)=\{[x,y]\colon x,y\in T\}$.

Since $\Gm/Z(\Gm)$ is non-cyclic (Lemma \ref{l:was-prop-F}) and has order $4$, it must be isomorphic to the Klein-four group, and hence has the the following properties:
\begin{itemize}
\item each non-identity element has order $2$;
\item any two elements commute with each other;
\item multiplying two distinct non-identity elements gives the third non-identity element.
\end{itemize}
Choose a transversal $\{e,a,b,c\}$ for $Z(\Gm)$ in $\Gm$.
Then the properties listed above imply:
\[
a\sim a^{-1}, b\sim b^{-1}, c\sim c^{-1} \quad\text{and}\quad
ab\sim c\sim ba, bc \sim a \sim cb, ca \sim b \sim ac.
\]
Since $a\sim a^{-1}$ and $b\sim ac$, we have
\[
[a,b] = [a^{-1},ac] = a^{-1}(ac) a (ac)^{-1} = [c,a];
\]
and since $a\sim cb\sim cb^{-1}$, we also have
\[
[a,b] = [cb^{-1},b] = (cb^{-1})b(bc^{-1})b^{-1}=[c,b].
\]
Thus $[a,b]=[c,a]=[c,b]$. By symmetry, these identities remain valid under any permutation of the symbols $a,b,c$, so we also obtain
$[b,c] = [a,b]=[a,c]$ and
$[c,a] = [b,c]=[b,a]$.
We have shown that $[a,b]$, $[b,c]$, $[c,a]$, $[b,a]$, $[c,b]$, $[a,c]$ are all equal; denoting their common value by~$z_0$, it follows from our previous remarks that
\[
\comm(\Gm) = \{ [x,y] \colon x,y\in T \} = \{e,z_0\}.
\]
Thus $|\comm(\Gm)|\leq 2$, and equality must hold since $\Gm$ is non-abelian.

Now let $\pi\in\Gmhat$ with $d_\pi\geq 2$; we must show that $H_\pi$ is $2$-dimensional.
Pick $x$ and $y$ in $\Gm$ which do not commute. Their images in $\Gm/Z(\Gm)$ are distinct, and so they generate all of $\Gm/Z(\Gm)$, by the properties of the Klein-four group mentioned above. Hence $\Gm$ is generated by $\{x,y\}\cup Z(\Gm)$.
Moreover, $\pi(Z(\Gm))\subseteq\Cplx I_\pi$ (by Schur's lemma), and so $\pi(\Gm)$ is generated by $\{ \pi(x), \pi(y)\}\cup\Cplx I_\pi$. Therefore, if we can produce a $2$-dimensional subspace $V$ which is $\pi(x)$-invariant and $\pi(y)$-invariant, this will force $V=H_\pi$ by irreducibility.

Let $v_1$ be an eigenvector of $\pi(y)$, with eigenvalue $\lambda$, and let $v_2=\pi(x)v_1$.
Since $\Gm/Z(\Gm)$ is commutative, $[x,y]$ is central in $\Gm$, and so by Schur's lemma $\pi([x,y])=\mu I_\pi$ for some\footnotemark\ $\mu\in \Cplx$.
\footnotetext{In fact, it follows from the proof of Lemma \ref{l:was-lem-C} that $\mu=-1$, but this is more than we need.}%
Thus $\pi(x)\pi(y)=\mu\pi(y)\pi(x)$, which implies $\pi(y)v_2=\lambda\mu v_1$.
Also, since every element of $\Gm/Z(\Gm)$ squares to the identity, $x^2\in Z(\Gm)$ and once again Schur's lemma yields $\pi(x)^2 \in \Cplx I_\pi$.
Hence $\pi(x)v_2 =\pi(x^2)v_1 \in \Cplx v_1$. Thus $V\defeq\operatorname{lin}\{v_1,v_2\}$ is an invariant subspace for $\pi(x)$ and $\pi(y)$, as required.
\end{proof}

\begin{proof}[Proof of Theorem \ref{t:was-G}\ref{li:was-G2}$\implies$\ref{li:was-G3}]
This follows immediately from the fact that $\Gmhat$ separates the points of~$\Gm$ (note that this does not require any form of the Plancherel theorem).
\end{proof}

To complete the circle of implications in Theorem \ref{t:was-G}, we need a closer look at $\pi(\Gm)$ when $\pi\in\Gmhat$ with $d_\pi=2$. For clarity, we isolate some steps in a preliminary lemma.
We denote by $Z_\Gm(x)$ the centralizer of an element $x\in\Gm$.

\begin{lem}\label{l:was-cor-D}
Let $\Gm$ be a group with $|\comm(\Gm)|=2$ and let $\pi\in\Gmhat$ with $d_\pi=2$. Let $x\in \Gm\setminus Z(\Gm)$.
\begin{romnum}
\item\label{li:was-D2}
If $c\in Z_\Gm(x)\setminus Z(\Gm)$, then $\pi(c)$ is a scalar multiple of~$\pi(x)$.
\item\label{li:was-D3}
If $a,b\in \Gm\setminus Z_\Gm(x)$, then $\pi(a^{-1}b)$ commutes with~$\pi(x)$.
\end{romnum}
\end{lem}

\begin{proof}
We start by diagonalizing the matrix $\pi(x)$.
Since $\pi(x)$ is unitary, its eigenvalues lie on the unit circle. On the other hand, pick $y\in \Gm$ with $xy\neq yx$.
Then by Lemma \ref{l:was-lem-C}\ref{li:was-C2}, $\pi(y^{-1})\pi(x)\pi(y)=-\pi(x)$; taking traces on both sides and rearranging yields $\Tr\pi(x)=0$. Therefore $\pi(x)$ has distinct eigenvalues $\lm$ and $-\lm$ for some $\lm\in\bbT$, and by taking corresponding eigenvectors $v_+$ and $v_{-}$ we obtain a basis for $H_\pi$.
If $c\in Z_\Gm(x)\setminus Z(\Gm)$, then since $\pi(c)$ commutes with $\pi(x)$ it must be diagonal with respect to the basis $\{v_+,v_{-}\}$. But since $c$ is non-central, repeating the argument used for $\pi(x)$ shows that the eigenvalues of $\pi(c)$ lie on the unit circle and sum to zero.
Thus $\pi(c)v_{\pm} = \pm \mu v_{\pm}$ for some $\mu\in\bbT$, so that $\pi(c) =\lm^{-1}\mu \pi(x)$, and we have proved part \ref{li:was-D2}.

For part \ref{li:was-D3}, it is convenient to introduce the matrix $J={\small\twomat{1}{0}{0}{-1}}$, so that $\pi(x)$ is represented by the matrix $\lm J$ with respect to the basis $\{v_+, v_-\}$. Now let $a,b\in \Gm$ with $ax\neq xa$ and $bx\neq xb$.
By Lemma \ref{l:was-lem-C}\ref{li:was-C2}, the operators $\pi(a)$ and $\pi(x)$ \emph{anti-commute}. A~direct calculation shows that an element of $M_2(\Cplx)$ which anti-commutes with $J$ must have zero entries on the diagonal, and so $\pi(a)$ swaps the two subspaces $\Cplx v_+$ and $\Cplx v_{-}$. The same is true for $\pi(b)$, and so $\pi(a^{-1}b)=\pi(a)^{-1}\pi(b)$ preserves each of these subspaces, i.e.\ it is diagonal with respect to the basis $\{v_+,v_{-}\}$. Hence $\pi(a^{-1}b)$ commutes with $\pi(x)$, and part \ref{li:was-D3} is proved.
\end{proof}

\begin{proof}[Proof of Theorem \ref{t:was-G}\ref{li:was-G3}$\implies$\ref{li:was-G1}]
Pick a non-central $x\in\Gm$. We will show that $\grindex{Z_\Gm(x)}{Z(\Gm)}=2$ and $\grindex{\Gm}{Z_\Gm(x)}=2$, which together imply that $\Zindex{\Gm}=4$.

For convenience let $\Om_n=\{\pi\in\Gmhat \colon d_\pi=n\}$.
We know that $Z_\Gm(x)\supseteq Z(\Gm)\cup\ xZ(\Gm)$. To prove the converse inclusion, let $c\in Z_\Gm(x)\setminus Z(\Gm)$. It suffices to prove that $x^{-1}c\in Z(\Gm)$, and since $\Om_1\cup\Om_2$ separates the points of $\Gm$, it suffices to prove that
$\pi(x^{-1}c)\in\Cplx I_\pi$ for all $\pi\in\Om_1\cup\Om_2$. This holds trivially if $\pi\in \Om_1$, and for $\pi\in\Om_2$ it follows from Lemma \ref{l:was-cor-D}\ref{li:was-D2}.

Similarly, to prove that $\grindex{\Gm}{Z_\Gm(x)}=2$, let $a$ and $b$ be elements of $\Gm\setminus Z_\Gm(x)$. It suffices to prove that
$\pi(a^{-1}bx)=\pi(xa^{-1}b)$ for all $\pi\in\Om_1\cup\Om_2$.
This holds trivially if $\pi\in\Om_1$, and for $\pi\in\Om_2$ it follows from Lemma \ref{l:was-cor-D}\ref{li:was-D3}.
\end{proof}

We can now give a proof of Theorem \ref{t:minimizers of AD} for countable groups.
\begin{itemize}
\item
Let $\Lm$ be a countable group satisfying $\AD(\Lm)=3/2$.
Since $\AD(\Lm)<\infty$, $\Lm$ is virtually abelian. Let $\nu$ be Plancherel measure for $\Lmhat$.
By Lemma \ref{l:was-A} $\nu(\Om_1)=1/2$ and $\nu(\{\pi\in\Lmhat\colon d_\pi \geq 3\})=0$.
Hence, by Proposition \ref{p:was-E} $|\comm(\Lm)|=2$. Moreover, if $x,x'\in\Lm$ satisfy $\pi(x)=\pi(x')$ for all $\pi\in\Lmhat$ with $d_\pi\leq 2$, then by injectivity of the non-abelian Fourier transform (or, alternatively, by applying the Plancherel identity \eqref{eq:plancherel-identity} to $\ip{\delta_x}{\delta_{x'}}$) we must have $x=x'$. By \ref{li:was-G3}$\implies$\ref{li:was-G2} in Theorem \ref{t:was-G}, $|\Lm:Z(\Lm)|=4$.

\item
Conversely, suppose $\Lm$ is a countable group with $|\Lm:Z(\Lm)|=4$. By \ref{li:was-G1}$\implies$\ref{li:was-G2} in Theorem \ref{t:was-G}, $|\comm(\Lm)|=2$ and $\maxdeg(\Lm)=2$. In particular, $\{\pi\colon d_\pi\geq 3\}$ is empty and so has measure zero, while by Proposition \ref{p:was-E} $\nu(\Om_1)=1/2$. Applying Lemma \ref{l:was-A} in the other direction we conclude that $\AD(\Lm)=3/2$. 
\end{itemize}

To obtain the general case we need another group-theoretic lemma.

\begin{lem}[Control by finitely generated subgroups]\label{l:Z}
Let $\Gm$ be a group and let $m\in\Nat$. If~$\Zindex{\Lm}\leq m$ for every finitely generated subgroup $\Lm\leq\Gm$, then $\Zindex{\Gm}\leq m$.
\end{lem}

\begin{proof}
We prove the contrapositive statement. Suppose $\Zindex{\Gm} \geq m+1$, and pick $a_1$,\dots, $a_{m+1}$ in $\Gm$ which belong to different cosets of $Z(\Gm)$. Then for each $(i,j)\in\Nat^2$ with $1\leq i  < j \leq m+1$, we have $a_i^{-1}a_j\notin Z(\Gm)$, so there exists $t_{ij}\in\Gm$ such that $a_i^{-1}a_j$ does not commute with~$t_{ij}$.
Let $\Lm$ be the subgroup of $\Gm$ generated by $\{a_1,\dots, a_{m+1}\}\cup\{t_{ij} \colon 1\leq i < j \leq m+1\}$. Then $\Lm$ is a finitely generated subgroup of $\Gm$, and $\Zindex{\Lm} \geq m+1$ since $a_1,\dots, a_{m+1}$ belong to different cosets of $Z(\Lm)$.
\end{proof}

\begin{proof}[Proof of  Theorem \ref{t:minimizers of AD} in the general case]
Suppose $\Gm$ is a group with $\AD(\Gm)=3/2$.
Since $\Gm$ is non-abelian, $\Zindex{\Gm}\geq 4$ by Lemma \ref{l:was-prop-F}. It suffices to show that $\Zindex{\Gm}\leq 4$.
Let $\Lm\leq \Gm$ be a finitely generated subgroup. If $\Lm$ is abelian $\Zindex{\Lm}=1$. If $\Lm$ is non-abelian, then by our sharp lower bound (Theorem \ref{t:sharp LB}) and monotonicity of $\AD$ (Proposition \ref{p:AD gen prop}), $3/2\leq \AD(\Lm)\leq\AD(\Gm)=3/2$.
So $\AD(\Lm)=3/2$, and since $\Lm$ is countable, the previous reasoning tells us that $\Zindex{\Lm}=4$.
Thus $\Zindex{\Lm}\leq 4$ for every finitely generated subgroup $\Lm\leq\Gm$. By Lemma \ref{l:Z} $\Zindex{\Gm}\leq 4$, as required.

Conversely, suppose $\Gm$ is a group satisfying $\Zindex{\Gm}=4$, and pick $x,y\in\Gm$ which do not commute. By countable saturation for $\AD$ (Proposition~\ref{p:AD countable-sat}) there is a countable subgroup $\Lm_0\leq\Gm$ such that $\AD(\Lm_0)=\AD(\Gm)$; taking $\Lm$ to be the subgroup generated by $\Lm_0$ and $\{x,y\}$, we see that $\Lm\leq\Gm$ is countable and non-abelian, and $\AD(\Lm)=\AD(\Gm)$ by monotonicity of~$\AD$.
Hence $\Zindex{\Lm}\geq 4$  by Lemma~\ref{l:was-prop-F}.
On the other hand, since $Z(\Lm)\supseteq \Lm\cap Z(\Gm)$,
\[
\Zindex{\Lm} \leq \grindex{\Lm}{\Lm\cap Z(\Gm)} \leq \Zindex{\Gm}  =4,
\]
and so $\Zindex{\Lm}=4$. Since $\Lm$ is countable, the previous reasoning tells us that $\AD(\Lm)=3/2$, and so $\AD(\Gm)=3/2$ as required.
\end{proof}

\end{section}


\begin{section}{Examples where $\AD(G)$ is as large as possible}
In this section we prove Theorem~\ref{t:all equal}. We shall build up to this result in stages.
Throughout this section, $[\Gm,\Gm]$ denotes the subgroup of $\Gm$ generated by the set of commutators (the \dt{commutator subgroup} or \dt{derived subgroup} of~$\Gm$).

\begin{lem}\label{l:ensures-char-null}
Let $\Gm$ be a countable VA group, equipped with counting measure, and let $\nu$ be the corresponding Plancherel measure on $\Gmhat$. Let $\Om_1=\{\pi\in\Gmhat \colon d_\pi =1 \}$.
If $[\Gm,\Gm]$ is infinite then $\nu(\Om_1)=0$.
\end{lem}


\begin{proof}
We prove the contrapositive statement. Assume $\Om_1$ has strictly positive measure, and for each $x\in\Gm$ put $h(x) = \int_{\Om_1} \pi(x)\,d\nu(\pi)$. Note that $h(x)=\nu(\Om_1)>0$ for all $x\in [\Gm,\Gm]$. We claim that $h\in \ell^2(\Gm)$. Assuming for the moment that this holds, we deduce that
\[
0 < | [\Gm,\Gm] |^{1/2}  \nu(\Om_1) \leq \norm{h}_2 < \infty
\]
and conclude that $| [\Gm,\Gm]| <\infty$ as required.

The claim follows from non-abelian Fourier inversion for $\Gm$, but we can give a proof that only requires\footnotemark\ Equation~\eqref{eq:plancherel-measure}. \footnotetext{That is, we only need the fact that the non-abelian Fourier transform is an isometry from $\ell^2(\Gm)$ to a suitable space of operator fields, rather than the subtler issue of characterizing the range of this isometry.}
Let $f\in c_{00}(\Gm)$, and note that by interchanging finite summation with integration we have
\[
\sum_{x\in\Gm} f(x) h(x)
= \sum_{x\in\Gm} f(x) \int_{\Om_1} \pi(x)\,d\nu(\pi)
= \int_{\Om_1} \pi(f) \,d\nu(\pi).
\]
Therefore, applying Cauchy--Schwarz and Equation~\eqref{eq:plancherel-measure} yields
\[
\left\vert \sum_{x\in \Gm} f(x) h(x) \right\vert
 \leq \nu(\Om_1)^{1/2} \left( \int_{\Om_1} |\pi(f)|^2\,d\nu(\pi) \right)^{1/2} 
 \leq \nu(\Om_1)^{1/2} \norm{f}_2 \;.
\]
By standard duality arguments in $\ell^2(\Gm)$ we deduce that $\norm{h}_2 \leq \nu(\Om_1)^{1/2} <\infty$, as claimed. This completes the proof.
\end{proof}

\begin{rem}
The proof of Lemma~\ref{l:ensures-char-null} can be generalized to any 2nd-countable unimodular Type~I group~$G$, and yields the following result: if the commutator subgroup in $G$ has non-compact closure, then the set of one-dimensional representations is a Plancherel-null subset of~$\Ghat$. We leave the details to the reader, since we have no applications at present for this more general result.
\end{rem}

\begin{prop}\label{p:AD=maxdeg}
Let $\Gm$ be a countable VA group. Suppose that the following two conditions hold:
\begin{romnum}
\item\label{li:one}
$[\Gm,\Gm]$ is infinite;
\item\label{li:two}
there exists some $d\geq 2$ such that $d_\pi \in\{1,d\}$ for all $\pi\in\Gmhat\}$.
\end{romnum}
Then $d=\AD(\Gm)=\AM(\FA(\Gm))=\maxdeg(\Gm)$.
\end{prop}

(Example \ref{eg:AD red-heis-mod2} shows that the assumption \ref{li:one} is necessary.)

\begin{proof}
Recall (see \eqref{eq:AD-AM-maxdeg} above) that by results in \cite{Run_PAMS06} we have $d=\maxdeg(\Gm)\geq \AM(\FA(\Gm)\geq \AD(\Gm)$.
It therefore suffices to prove that $\AD(\Gm)=d$.

Let $\nu$ be the Plancherel measure on $\Gmhat$ corresponding to counting measure on $\Gm$. By Assumption~\ref{li:two},
\[
1 = \int_{\Gmhat} d_\pi \,d\nu(\pi) = \nu(\Om_1)+d\nu(\Om_d),
\]
while by Theorem~\ref{t:main formula},
\[
\AD(\Gm) = \int_{\Gmhat} d_\pi^2 \,d\nu(\pi) = \nu(\Om_1)+d^2\nu(\Om_d).
\]
Hence
$\AD(\Gm) = \nu(\Om_1) + d(1-\nu(\Om_1))$. But, by Assumption~\ref{li:one} and Lemma~\ref{l:ensures-char-null}, $\nu(\Om_1)=0$, and the result follows.
\end{proof}

\begin{lem}\label{l:nice semidirect}
Let $N$ be a torsion-free LCA group, let $m\geq 2$, and let $\alpha:  C_m\to {\rm Aut}(N)$. Let $G=N\rtimes_\alpha  C_m$ be the corresponding semidirect product. Then $d_\pi \mid m$ for all $\pi\in\Ghat$.
\end{lem}

\begin{proof}
This follows from (a simple case of) the Mackey machine, as described in e.g.~\cite{KT_book}.
Given $\chi\in\widehat{N}$, let $D_\chi$ denote the stabilizer of $\chi$ with respect to the adjoint action of $C_p$ on~$\widehat{N}$, and let $G_\chi = N\rtimes_\alpha D_\chi \subset G$. Then each $\pi\in\Ghat$ is unitarily equivalent to one obtained by inducing some one-dimensional representation $G_\chi\to\bbT$. (See \cite[Theorem 4.40]{KT_book} for a more general result.) In particular,
$d_\pi = \grindex{G}{G_\chi} = \grindex{C_m}{D_\chi}$, which divides $| C_m|=m$ by Lagrange's theorem.
\end{proof}

\begin{proof}[The proof of Theorem \ref{t:all equal}]
Recall the setup: $N$ is a torsion-free LCA group, $p$ is a prime, and $G=N\rtimes_\alpha C_p$ for some non-trivial action~$\alpha$.
By Lemma \ref{l:nice semidirect} we have $\maxdeg(G)\leq p$. (In fact, equality holds, but we will obtain this for free later on.)

By assumption there exists $x_0\in N$ whose $\alpha$-orbit contains at least one other point of $N$. Let $N_{00}$ be the subgroup of $N$ (not necessarily closed!) that is generated by this orbit; then $N_{00}$ is countable, torsion-free abelian, and $\alpha$-invariant.
Let $\Gm=N_{00}\rtimes_\alpha C_p$: this is a subgroup of $G_d$, and so by monotonicity of $\AD$,
\begin{equation}\label{eq:chain-will-be-equal}
\AD(\Gm) \leq \AD(G_d) \leq \AM(\FA(G)) \leq p.
\tag{$*$}
\end{equation}

Now observe the following:
\begin{romnum}
\item
Let $x_1$ be any point in the $C_p$-orbit of $x_0$ with $x_1\neq x_0$. Since  $x_0$ and $x_1$ are conjugate in $\Gm$, $x_0^{-1}x_1$ is a non-trivial commutator in $\Gm$; also $x_0^{-1}x_1\in N_{00}$, which is torsion-free. Thus $[\Gm,\Gm]$ contains a free abelian subgroup, and in particular is infinite. 
\item
Applying Lemma \ref{l:nice semidirect} with $N$ replaced by $N_{00}$, and recalling that $p$ is prime, we see that $d_\pi \in\{1,p\}$ for every $\pi\in\Gmhat$.
\end{romnum}
Applying Proposition~\ref{p:AD=maxdeg} we conclude that $\AD(\Gm)=p$, and so all the inequalities in \eqref{eq:chain-will-be-equal} must be equalities.
\end{proof}

\end{section}

%

 \begin{section}{Concluding remarks and further questions}
The most obvious and important question is:

\begin{qn}\label{q:1}
Do we always have $\AM(\FA(G))=\AD(G_d)$?
\end{qn}

In light of the good functorial properties of $\AD$, a positive answer would provide new tools for estimating or calculating amenability constants of Fourier algebras (c.f.~Remark \ref{r:AD-better-than-AMA}).

To see why the answer to Question \ref{q:1} is not obvious, we give a sketch of how the inequality $\AM(\FA(G))\geq \AD(G_d)$ is proved in \cite{FR_MZ05, Run_PAMS06}.
By definition, $\AM(\FA(G))$ is the infimum of norms of virtual diagonals for $\FA(G)$. We will not define virtual diagonals here: they are elements ${\bf M}\in (\FA(G)\ptp \FA(G))^{**}$ satisfying certain algebraic properties, where $\ptp$ denotes the projective tensor product of Banach spaces. There is a linear contraction
\[
J: \FA(G)\ptp\FA(G)\to \FA(G\times G) \;, \quad J(f\tp g)(x,y) = f(x)g(y^{-1})\,,
\]
which extends to a linear contraction $\widetilde{J}: (\FA(G)\ptp\FA(G))^{**} \to \FS(G_d\times G_d)$. Moreover, the algebraic properties of ${\bf M}$ imply that $\widetilde{J}({\bf M})= 1_{\adiag(G)}$. Hence
\[ \AD(G_d) = \fsnorm{1_{\adiag(G)}} \leq \inf_{\bf M} \norm{\bf M} = \AM(\FA(G)). \]
The point is that for non-abelian $G$, $J$ is never an isometry (although it is bijective if $G$ is virtually abelian), and $\widetilde{J}$ is not a quotient map of Banach spaces. Hence even if one knows that the set $\widetilde{J}^{-1}(\{ 1_{\adiag(G)}\})$ contains a virtual diagonal ${\bf M}$, there is no immediate reason why one should be able to find such an ${\bf M}$ whose norm equals that of $1_{\adiag(G)}$.

On the other hand, if $G$ is finite and non-abelian, then even though $J:\FA(G)\ptp\FA(G)\to\FA(G\times G)$ need not be isometric we know from Theorem \ref{t:AD of finite} that we do have $\AM(\FA(G))=\AD(G)$. So if there is a group $G$ for which $\AM(\FA(G))>\AD(G_d)$, there must be some deeper obstruction behind this which is not visible for finite groups; Theorem \ref{t:all equal} shows that we can have groups that look very different from (abelian)$\times$(finite) where $\AM(\FA(G))=\AD(G_d)$.

An additional motivation for Question \ref{q:1} is that a positive answer, combined with Theorem~\ref{t:minimizers of AD}, would give a complete characterization of those locally compact non-abelian $G$ where $\AM(\FA(G))$ attains its minimal value. We phrase this as a separate question, which might be easier to attack directly.

\begin{qn}\label{q:2}
Which non-abelian $G$ attain the minimal value for $\AM(\FA(G))$?
\end{qn}

If $\AM(\FA(G))=3/2$ then, since $G$ is non-abelian, our sharp lower bound (Theorem \ref{t:sharp LB}) implies that
\[
3/2\leq \AD(G_d)\leq \AM(\FA(G)) =3/2
\]
and so by Theorem \ref{t:minimizers of AD} $\Zindex{G}=4$. Thus Question \ref{q:2} is really asking if the converse holds, i.e.\ does $\Zindex{G}=4$ imply $\AM(\FA(G))=3/2$?
(For finite groups, the answer is positive by Theorem \ref{t:AD of finite}.)

\begin{qn}\label{q:3}
Can we classify those $\Gm$ such that $3/2 < \AD(\Gm) < 2$?
\end{qn}

Note that  if $\Lm$ is countable and $|[\Lm,\Lm]|$ is infinite, then by Lemma \ref{l:ensures-char-null} and the main formula (Equation \eqref{eq:AD-explicit-formula}) we have $\AD(\Lm)\geq 2$. Hence the condition that $\AD(\Gm)<2$ forces $[\Gm,\Gm]$ to be finite (if it were infinite then there would be a countable $\Lm\leq\Gm$ witnessing this). There is a substantial literature on groups whose commutator subgroups are finite (so-called FD groups) and one might be able to use known structure theory to attack Question~\ref{q:3}.
Progress might also result from a good answer to the following problem.

\begin{qn}\label{q:4}
Can we bound $\maxdeg(G)$ from above by some explicit function of $\AD(G_d)$?
\end{qn}

This was done for finite groups in \cite[Theorem 4.9]{LLW_SM96} but to the author's knowledge it remains an open problem for infinite VANA groups.

\end{section}

\section*{Acknowledgements}
The challenge of finding optimal lower bounds on $\AM(\FA(G))$ for infinite non-abelian groups originally arose during conversations in 2016 with M.~Alaghmandan, who the author thanks for some stimulating discussions.
He also thanks M.~Ghandehari for several valuable discussions about the overall project, and is indebted to her for the crucial suggestion to make use of results in \cite[Section~3G]{Ars_Lyon76}.
Thanks also go to: Y.~Cornulier for pointing out via {\itshape MathOverflow} that an example which the author had laboriously constructed by hand was in fact isomorphic to the group in Example~\ref{eg:red-heis-mod2}; N.~Juselius for sharing a copy of~\cite{Juselius_MSc}; and K.~F. Taylor for sharing some private notes on the Plancherel theorem for crystal groups, which led to the additional examples mentioned in Remark \ref{r:KFT_crystal}.

Key progress towards the main formula (Theorem~\ref{t:main formula}) was made, and some special cases of Theorem~\ref{t:all equal} were obtained, through discussions with Ghandehari during a visit to the University of Delaware in March 2020. The author thanks his hosts for their hospitality, and acknowledges financial support from the Visitor Fund of the Department of Mathematics and Statistics at Lancaster University. Important additions and revisions were worked out while the author was visiting the Laboratoire de math\'ematiques de Besan\c{c}on, Universit\'e de Franche-Comt\'e, in November 2021. The author thanks his hosts for a welcoming and pleasant working environment, and acknowledges financial support from the Alliance Hubert Curien programme, as well as valuable moral support from N.~Blitvi\'c and J.~M. Lindsay.

Finally, thanks are due to an anonymous referee for their careful reading, and for several helpful and detailed suggestions, which led to the replacement of some rather artificial {\it ad hoc} arguments with the general perspective described in Section \ref{ss:dual-of-product} and the Appendix.

\appendix

\begin{section}{Comparing $\hatonehattwo{G}$ and $\hatonetwo{G}$ as topological spaces}\label{app:general-dual-of-product}
The main purpose of this appendix is to show how the result for general locally compact groups that was stated in Remark \ref{r:proposed by referee} can be derived from standard ingredients in \cite{DixCstar}.
(However, for acknowledgement of other coverage in the literature, see the author's note at the end of Section \ref{ss:dual-of-product}.)
We restate the desired result for the reader's convenience.

\begin{prop}\label{p:RTP_group}
Let $G_1$ and $G_2$ be locally compact groups, and let $J:\hatonehattwo{G}\to \hatonetwo{G}$ be defined by $J([\pi_1],[\pi_2]) = [\pi_1\tp \pi_2]$.  Then $J$ is continuous and relatively open.
\end{prop}

This result was suggested to the author by the referee, who observed that it is a special case of general results concerning states on $\Cst$-algebras and the corresponding GNS representations. To fix notation: for a $\Cst$-algebra~$A$, we write $S(A)$ for the set of states on~$A$, equipped with the relative \wstar-topology inherited as a subset of~$A^*$. Given $\Cst$-algebras $A$ and $B$ and states $f\in S(A)$ and $g\in S(B)$, the linear functional $f\tp g :A \tp B \to \Cplx$ extends uniquely to a state on the maximal $\Cst$-tensor product $A\maxtp B$. We then have the following result.

\begin{prop}\label{p:compare state spaces}
Let $A$ and $B$ be $\Cst$-algebras and let $D= A \maxtp B$. The function $T:S(A)\times S(B) \to S(D)$ which sends $(f,g)$ to $f\tp g$ is a homeomorphism onto its range (where the range is given the subspace topology).
\end{prop}

The following proof is based on a sketch provided by the referee.
\begin{proof}
Let $((f_i,g_i))_{i\in\Lambda}$ be a net in $S(A)\times S(B)$ and let $(f,g)\in S(A)\times S(B)$.

Suppose $(f_i,g_i) \to (f,g)$ in $S(A)\times S(B)$. Then $f_i(a)\to f(a)$ and $g_i(b)\to g(b)$ for each $a\in A$ and $b\in B$, and so $(f_i\tp g_i)(w)\to (f\tp g)(w)$ for each $w\in A\tp B$. Since $A\tp B$ is norm-dense in $D$ and the net $(f_i\tp g_i)$ is bounded in $D^*$, it follows that $f_i\tp g_i \to f\tp g$ in $S(D)$. This shows that $T$ is continuous.

We now claim that the converse holds: if $f_i\tp g_i \to f\tp g$ in $S(D)$, then $f_i\to f$ in $S(A)$ and $g_i\to g$ in $S(B)$. Assuming this claim for the moment: by considering nets that only take two values, we see that $T$ is injective (since we are dealing with Hausdorff spaces); and then by considering general nets, we see that $T^{-1} : T(S(A)\times S(B)) \to S(A)\times S(B)$ is continuous, as required.

It remains to prove the claim. Let $a\in A$ and $b\in B$: we aim to show that $f_i(a)\to f (a)$ and $g_i(b)\to g(b)$. To motivate the proofs, note that the special case where $A$ and $B$ are unital follows immediately from the following calculations:
\[ \begin{aligned}
f_i(a) & = (f_i\tp g_i)(a\tp 1_B) \to (f\tp g)(a\tp 1_B) = f(a), \\
g_i(b) & = (f_i\tp g_i)(1_A\tp b) \to (f\tp g)(1_A\tp b) = g(b).
\end{aligned} \]

For the general case, let $A^+_1$ be the positive part of the closed unit ball of $A$, and define $B^+_1$ similarly. Given $\veps \in (0,1/4)$, choose $x\in A^+_1$ and $y\in B^+_1$ such that $f(x), g(y) \in [1-\veps,1]$. Since $f\tp g= \wstar\lim_i f_i\tp g_i$, we have
\[ \begin{aligned}
1\geq f_i(x)g_i(y)
& \geq f(x)g(y)-\veps^2 \\
& \geq (1-\veps)^2-\veps^2 = 1-2\veps
\quad\text{for all sufficiently large $i$.}
\end{aligned} \]
Since $f_i(x)g_i(y) \leq g_i(y) \leq 1$, it follows that $g_i(y) \in [1-2\veps,1]$ for all sufficiently large~$i$.
Hence, for all such $i$,
\[ \begin{aligned}
|f_i(a)-f(a)|
& \leq |f_i(a)-f_i(a)g_i(y)| + |f_i(a)g_i(y)-f(a)g(y) | + |f(a)g(y)-f(a)| \\
& \leq  2\veps \norm{a}   + |f_i(a)g_i(y)-f(a)g(y) | +  \veps\norm{a} \\
\end{aligned} \]
where the middle term tends to $0$. Thus $|f_i(a)-f(a)| \leq (3\norm{a} + 1)\veps$ for all sufficiently large~$i$. This proves that $f(a)=\lim_i f_i(a)$, and a similar argument shows that $g(b)=\lim_i g_i(b)$.
\end{proof}

\begin{rem}
The argument above works if $D$ is replaced by the completion of $A\tp B$ in any $\Cst$-tensor norm, so it might seem more natural to use the minimal tensor product. Our choice is motivated by the application to (full) group $\Cst$-algebras: if $G_1$ and $G_2$ are locally compact groups, then $\Cst(G_1\times G_2) \cong \Cst(G_1)\maxtp\Cst(G_2)$.
\end{rem}

To avoid ambiguity in what follows: a~``representation'' of a $\Cst$-algebra $A$ means a non-degenerate $*$-homomorphism $A\to\Bdd(H)$ for some Hilbert space; ``irreducible'' means topologically irreducible; and ``equivalence'' of representations means unitary equivalence. We then write $\Ahat$ for the set of equivalence classes of irreducible representations of $A$; this carries a canonical topology.

In the remainder of this appendix, we show how Proposition \ref{p:compare state spaces} implies a corresponding result for the topological spaces $\Ahat$, $\Bhat$ and~$\Dhat$ (Corollary \ref{c:RTP_Cstar} below), of which Proposition \ref{p:RTP_group} is a special case.
We start with a quick sketch of basic facts concerning tensor products of representations, whose proofs we could not find in  \cite{DixCstar}.
Given representations $\pi:A\to \Bdd(H)$ and $\sigma: B\to\Bdd(K)$, the homomorphism $\pi\tp\sigma: A\tp B \to \Bdd(H\tp_2 K)$ extends uniquely to a representation of $D=A\maxtp B$ (other $\Cst$-tensor norms work just as well). Importantly, if $\pi$ and $\sigma$ are irreducible then $\pi\tp\sigma$ is also irreducible. The easiest way to see this, given the material in \cite{DixCstar}, is as follows:
\begin{itemize}
\item by a version of Schur's lemma for representations of $\Cst$-algebras (see \cite[Proposition 2.3.1]{DixCstar}, we have $\pi(A)'=\Cplx I_H$ and $\sigma(B)'=\Cplx I_K$;
\item since $\pi$ and $\sigma$ are non-degenerate representations, an approximation argument shows that the SOT-closure of $\pi(A)\tp \sigma(B)$ inside $\Bdd(H\tp_2K)$ contains both $\pi(A)\tp\Cplx I_K$ and $\Cplx I_H \tp \sigma(A)$;
\item by basic calculations with commutants, we deduce that $(\pi\tp \sigma)(D)' = \Cplx {I}_{H\tp_2 K}$, and then apply the other direction of Schur's lemma.
\end{itemize}
We therefore obtain a well-defined function
\begin{equation}\label{eq:J-for-Cstar}
J:\Ahat\times\Bhat\to \Dhat \;, \quad ([\pi],[\sigma]) \mapsto [\pi\tp\sigma].
\end{equation}

Given $f\in S(A)$, the corresponding GNS representation $\pi_f: A \to \Bdd(H_f)$ is irreducible if and only if the state $f$ is pure  (\cite[Proposition 2.5.4]{DixCstar}). Therefore, writing $P(A)$ for the set of pure states on $A$, there is a well-defined function $q_A:P(A)\to \Ahat$, $f\mapsto [\pi_f]$. Furthermore, if $f\in P(A)$ and $g\in P(B)$, and $f\tp g$ denotes the corresponding state on $D$, then $\pi_{f\tp g}$ can be identified with $\pi_f\tp \pi_g : D \to \Bdd(H_f\tp_2 H_g)$ which is irreducible by our earlier remarks, and so $f\tp g$ is pure. Thus $T: S(A)\times S(B) \to S(D)$ restricts to a well-defined function $P(A)\times P(B) \to P(D)$, and we have a commutative diagram:

\begin{equation}\label{eq:CD}
\begin{tikzcd}
	{P(A)\times P(B)} && {P(D)} \\
	\\
	{\Ahat\times\Bhat} && {\Dhat}
	\arrow["{q_A\times q_B}"', from=1-1, to=3-1]
	\arrow["T", from=1-1, to=1-3]
	\arrow["J"', from=3-1, to=3-3]
	\arrow["{q_D}", from=1-3, to=3-3]
\end{tikzcd}
\end{equation}

Note that $q_A$ is surjective (\cite[Section 2.4.6]{DixCstar}). Moreover, when $P(A)$ is given the relative \wstar-topology and $\Ahat$ is given its canonical topology, $q_A$ is continuous and open (\cite[Theorem 3.4.11]{DixCstar}). Similarly $q_B$, $q_D$ and $q_A\times q_B$ are all continuous open surjections. We are therefore in a position to apply the following lemma.

\begin{lem}\label{l:gen-top-square}\
Let $X$, $X'$, $Y$ and $Y'$ be topological spaces. Suppose there exist continuous open surjections $q_X: X' \to X$ and $q_Y:Y'\to Y$, and suppose there exist functions $f':X'\to Y'$ and $f:X\to Y$ satisfying $q_Y f' = fq_X$.
\begin{romnum}
\item If $f'$ is continuous, then so is $f$.
\item If $f'$ is relatively open, then so is $f$.
\end{romnum}
\end{lem}

The proof is an elementary exercise in general point-set topology (no need for nets) and we leave it to the reader to check the details. Applying this lemma with $X'=P(A)\times P(B)$, $Y=\Dhat$, etc.\ and
appealing to Proposition \ref{p:compare state spaces}, we obtain the following result.

\begin{cor}\label{c:RTP_Cstar}
The function $J$ from Equation \eqref{eq:J-for-Cstar} is continuous and relatively open.
\end{cor}

Finally, let us specialise to the case where $A=\Cst(G_1)$ and $B=\Cst(G_2)$ for arbitrary locally compact groups $G_1$ and $G_2$, and $D=\Cst(G_1\times G_2) = \Cst(G_1)\maxtp \Cst(G_2)$.
As explained in \cite[Section 13.3]{DixCstar}, for any $G$ there is a natural identification of $\Ghat$ with $\longhat{\Cst(G)}$, and  \emph{by definition}  the Fell topology on $\Ghat$ corresponds to the canonical topology on $\longhat{\Cst(G)}$. Thus Proposition \ref{p:RTP_group} is a special case of Corollary \ref{c:RTP_Cstar}.

\end{section}




\begin{thebibliography}{LLW96}

\bibitem[Ars76]{Ars_Lyon76}
Gilbert Arsac.
\newblock Sur l'espace de {B}anach engendr\'e par les coefficients d'une
  repr\'esentation unitaire.
\newblock {\em Publ. D\'ep. Math. (Lyon)}, 13(2):1--101, 1976.

\bibitem[Dix77]{DixCstar}
Jacques Dixmier.
\newblock {\em {$C\sp*$}-algebras}, volume~15 of {\em North-Holland
  Mathematical Library}.
\newblock North-Holland Publishing Co., Amsterdam, 1977.
\newblock Translated from the French by Francis Jellett.

\bibitem[Dix81]{DixVNA}
Jacques Dixmier.
\newblock {\em {V}on {N}eumann algebras}, volume~27 of {\em North-Holland
  Mathematical Library}.
\newblock North-Holland Publishing Co., Amsterdam, 1981.
\newblock Translated from the second French edition by F. Jellett.

\bibitem[Eym64]{eymard}
Pierre Eymard.
\newblock L'alg\`ebre de {F}ourier d'un groupe localement compact.
\newblock {\em Bull. Soc. Math. France}, 92:181--236, 1964.

\bibitem[Fel63]{Fell63_wc-tp}
James M.~G. Fell.
\newblock Weak containment and {K}ronecker products of group representations.
\newblock {\em Pacific J. Math.}, 13:503--510, 1963.

\bibitem[Fol16]{Folland_ed2}
Gerald~B. Folland.
\newblock {\em A course in abstract harmonic analysis}.
\newblock Textbooks in Mathematics. CRC Press, Boca Raton, FL, second edition,
  2016.

\bibitem[FR05]{FR_MZ05}
Brian~E. Forrest and Volker Runde.
\newblock Amenability and weak amenability of the {F}ourier algebra.
\newblock {\em Math. Z.}, 250(4):731--744, 2005.

\bibitem[FR11]{FR_CMB11}
Brian~E. Forrest and Volker Runde.
\newblock Norm one idempotent {$cb$}-multipliers with applications to the
  {F}ourier algebra in the {$cb$}-multiplier norm.
\newblock {\em Canad. Math. Bull.}, 54(4):654--662, 2011.

\bibitem[Joh94]{BEJ_AG}
Barry~E. Johnson.
\newblock Non-amenability of the {F}ourier algebra of a compact group.
\newblock {\em J. London Math. Soc. (2)}, 50(2):361--374, 1994.

\bibitem[Jus19]{Juselius_MSc}
Nicholas~A. Juselius.
\newblock Bounds on the amenability constants of {F}ourier algebras in the
  completely bounded multiplier norm.
\newblock Master's thesis, University of Alberta, 2019.

\bibitem[KL18]{KL_AGBGbook}
Eberhard Kaniuth and Anthony To-Ming Lau.
\newblock {\em Fourier and {F}ourier-{S}tieltjes algebras on locally compact
  groups}, volume 231 of {\em Mathematical Surveys and Monographs}.
\newblock American Mathematical Society, Providence, RI, 2018.

\bibitem[KT13]{KT_book}
Eberhard Kaniuth and Keith~F. Taylor.
\newblock {\em Induced representations of locally compact groups}, volume 197
  of {\em Cambridge Tracts in Mathematics}.
\newblock Cambridge University Press, Cambridge, 2013.

\bibitem[LLW96]{LLW_SM96}
Anthony To-Ming Lau, Richard~J. Loy, and George~A. Willis.
\newblock Amenability of {B}anach and {$C^*$}-algebras on locally compact
  groups.
\newblock {\em Studia Math.}, 119(2):161--178, 1996.

\bibitem[MP16]{MP_JFA16}
Jayden Mudge and Hung~Le Pham.
\newblock Idempotents with small norms.
\newblock {\em J. Funct. Anal.}, 270(12):4597--4603, 2016.

\bibitem[Run06]{Run_PAMS06}
Volker Runde.
\newblock The amenability constant of the {F}ourier algebra.
\newblock {\em Proc. Amer. Math. Soc.}, 134(5):1473--1481, 2006.

\bibitem[RW98]{RW_morita-book}
Iain Raeburn and Dana~P. Williams.
\newblock {\em Morita equivalence and continuous-trace {$C^*$}-algebras},
  volume~60 of {\em Mathematical Surveys and Monographs}.
\newblock American Mathematical Society, Providence, RI, 1998.

\bibitem[Tak02]{Takesaki_vol1}
Masamichi Takesaki.
\newblock {\em Theory of operator algebras. {I}}, volume 124 of {\em
  Encyclopaedia of Mathematical Sciences}.
\newblock Springer-Verlag, Berlin, 2002.

\bibitem[Tho68]{Thoma68}
Elmar Thoma.
\newblock Eine {C}harakterisierung diskreter {G}ruppen vom {T}yp~{I}.
\newblock {\em Invent. Math.}, 6:190--196, 1968.

\bibitem[Wal72]{Wal_JFA72}
Martin~E. Walter.
\newblock {$W^{\ast} $}-algebras and nonabelian harmonic analysis.
\newblock {\em J. Functional Analysis}, 11:17--38, 1972.

\end{thebibliography}

\vfill

\newcommand{\address}[1]{{\small\sc#1.}}
\newcommand{\email}[1]{\texttt{#1}}

\noindent
\address{Yemon Choi,
Department of Mathematics and Statistics,
Lancaster University,
Lancaster LA1 4YF, United Kingdom} 

\noindent
\email{y.choi1@lancaster.ac.uk}

\end{document}